\documentclass[a4paper,11pt]{article}
\usepackage{graphicx}
\usepackage{geometry}
\usepackage{enumerate}
\usepackage{amsmath}
\usepackage{amsthm}
\usepackage{bbm}
\usepackage{bm}
\usepackage{amssymb}
\usepackage{mathrsfs}
\usepackage{dsfont}
\usepackage{float}
\usepackage{booktabs,multirow}
\usepackage{makecell}
\usepackage{setspace}
\usepackage{subfigure}
\usepackage[sort&compress, comma]{natbib}
\makeatletter
\renewcommand\@biblabel[1]{${#1}.$}
\usepackage[colorlinks=true,linkcolor=blue, citecolor=blue, urlcolor=blue]{hyperref}
\hypersetup{hypertex=true,
	colorlinks=true,
	linkcolor=blue,
	anchorcolor=blue,
	citecolor=blue}
\makeatother
\usepackage{color}
\usepackage[font=small,format=plain,labelfont=bf,up,textfont=normal,up,justification=justified,singlelinecheck=false]{caption}
\usepackage[normalem]{ulem}

\usepackage{pgfplots}

\usepackage{authblk}
\usepackage{tikz}
\usetikzlibrary{shapes.geometric, arrows}

\DeclareMathOperator*{\argmin}{argmin}

\usepackage[noend]{algpseudocode}
\usepackage{algorithmicx,algorithm}
\usepackage{epstopdf}
\usepackage{multirow}

\newtheorem{theorem}{Theorem}[section]
\newtheorem{corollary}{Corollary}[section]

\newtheorem{definition}{Definition}
\newtheorem{proposition}{Proposition}
\newtheorem{assumption}{Assumption}
\newtheorem{remark}{Remark}

\usepackage{ulem}

\begin{document}

\title{\textbf{High-Dimensional Binary Variates: Maximum Likelihood Estimation with Nonstationary Covariates and Factors}}

\author[1]{\normalsize Xinbing Kong}
\author[2]{\normalsize Bin Wu\thanks{Corresponding author. Email: bin.w@ustc.edu.cn}}
\author[2]{\normalsize Wuyi Ye}
\affil[1]{Southeast University, Nanjing 211189, China}
\affil[2]{University of Science and Technology of China, Hefei 230026, China}
\date{}

\maketitle

\begin{abstract}
This paper introduces a high-dimensional binary variate model that accommodates nonstationary covariates and factors, and studies their asymptotic theory. This framework encompasses scenarios where single indices are nonstationary or cointegrated. For nonstationary single indices, the maximum likelihood estimator (MLE) of the coefficients has dual convergence rates and is collectively consistent under the condition $T^{1/2}/N\to0$, as both the cross-sectional dimension $N$ and the time horizon $T$ approach infinity. The MLE of all nonstationary factors is consistent when $T^{\delta}/N\to0$, where $\delta$ depends on the link function. The limiting distributions of the factors depend on time $t$, governed by the convergence of the Hessian matrix to zero. In the case of cointegrated single indices, the MLEs of both factors and coefficients converge at a higher rate of $\min(\sqrt{N},\sqrt{T})$. A distinct feature compared to nonstationary single indices is that the dual rate of convergence of the coefficients increases from $(T^{1/4},T^{3/4})$ to $(T^{1/2},T)$. Moreover, the limiting distributions of the factors do not depend on $t$ in the cointegrated case. Monte Carlo simulations verify the accuracy of the estimates. In an empirical application, we analyze jump arrivals in financial markets using this model, extract jump arrival factors, and demonstrate their efficacy in large-cross-section asset pricing.
\end{abstract}

\newpage

\section{Introduction}\label{sec:introduction}
Since \cite{chamberlain1983arbitrage}'s introduction of the linear approximate factor model, factor models have been the focus of extensive research. Seminal contributions include those by \cite{bai2003inferential}, \cite{fan2013large}, \cite{pelger2019large}, and \cite{he2025huber}. In recent years, attention has increasingly turned toward nonlinear factor models. For example, \cite{chen2021quantile} developed estimation methods and asymptotic theory for quantile factor models. In particular, considerable works have addressed nonlinear binary factor models, c.f., \cite{chen2021nonlinear}, \cite{ando2022bayesian}, \cite{wang2022maximum}, \cite{ma2023generalized}, and \cite{gao2023binary}. Binary factor models have diverse applications in fields such as engineering (data compression, visualization, pattern recognition, and machine learning), economics and finance (credit default analysis, macroeconomic forecasting, jump arrival analysis), and biology (gene sequence analysis). However, the aforementioned studies assume that both the factors and covariates are stationary processes---an assumption that may not hold in practice. For instance, daily jump frequencies in financial markets often exhibit aggregation in jump occurrence probabilities (see \citealt{bollerslev2011estimation,bollerslev2011tails}), suggesting nonstationarity, which is also justified in our empirical studies. This paper addresses this gap by investigating a general class of binary factor models where both covariates and factors are generated by integrated processes.

Specifically, consider the following single-index general factor model:
\begin{align}\label{eq:single-index general factor model} y_{it}=\Psi(z_{0it})+u_{it}~~\text{where}~~z_{0it}=\beta_{0i}'x_{it}+\lambda_{0i}'f_{0t}, \ \mbox{for} \ i=1,..., N \ \mbox{and} \ t=1,..., T.
\end{align}
Here, $x_{it}$ is a $q$-dimensional covariate with coefficient vector $\beta_{0i}$, and $f_{0t}$ is an $r$-dimensional factor with factor loading vector $\lambda_{0i}$. The binary outcome $y_{it}$ is modeled through a known nonlinear link function $\Psi(\cdot)$ (such as logit or probit). Both $x_{it}$ and $f_{0t}$ are integrated of order one, i.e., $I(1)$ processes. It is certainly that $y_{it}$ can be extended to other types of variate like counts. We consider two cases for the single index $z_{0it}$, one is nonstationary $I(1)$ index and one is cointegrated $I(0)$ index. In the nonstationary univariate regression setting (i.e., $\lambda_{0i}'f_{0t}=0$ and $N=1$), \cite{park2000nonstationary} provide the relevant asymptotic theory. However, in high-dimensional settings with a factor structure, the asymptotic properties remain unexplored. Furthermore, when $z_{0it}$ is cointegrated, the associated asymptotic theory is also missing.

This paper develops a theoretical framework for binary factor models with nonstationary covariates and factors, modeled as integrated time series. Specifically, we model the covariates and factors as integrated time series and consider scenarios where the single index is either nonstationary or cointegrated. Our approach builds on earlier work on the asymptotics of nonlinear functions of integrated time series (c.f., \citealt{park1999asymptotics,park2000nonstationary,park2001nonlinear}; \citealt{dong2016estimation}; \citealt{zhou2024semi}) and on MLE methodologies for high-dimensional stationary factor models (c.f., \citealt{bai2012statistical}; \citealt{chen2021quantile}; \citealt{gao2023binary}; \citealt{yuan2023two}; \citealt{xu2025quasi}). There are also studies involving nonstationary time series in high-dimensional linear models, such as \cite{zhang2018clt}, \cite{dong2021varying}, \cite{trapani2021inferential}, and \cite{barigozzi2024inference}. Our main theoretical contribution lies in establishing the asymptotic properties of the MLEs for both coefficients and factors within the class of general binary factor models.

The findings of this paper are summarized below. In the model \eqref{eq:single-index general factor model} with a nonstationary single-index, the convergence rate of the estimated $\alpha_{0i}$ ($\alpha_{0i}=(\beta_{0i}',\lambda_{0i}')'$) is characterized by two distinct rates along different axes in a new coordinate system where $\alpha_{0i}/\|\alpha_{0i}\|$ defines one axis. Along the
axis of $\alpha_{0i}/\|\alpha_{0i}\|$, the estimated $\alpha_{0i}$ (denoted by $\hat{\alpha}_i$) converges at a rate of $T^{1/4}$, while along axes orthogonal to $\alpha_{0i}/\|\alpha_{0i}\|$, it converges at a faster rate of $T^{3/4}$. This dual-rate convergence is consistent with findings in binary regression models, such as the univariate case documented in \cite{park2000nonstationary}. Collectively, $\hat{\alpha}_{i}$ exhibits a convergence rate of $T^{1/4}$. The convergence rate of $\hat{f}_t$ to $f_{0t}$ is $t^{-\delta/2}N^{1/2}$, where $\delta$ is determined by the property of the link function. The normalized estimator $\hat{\alpha}_{i}/\|\hat{\alpha}_{i}\|$ converges to $\alpha_{0i}/\|\alpha_{0i}\|$ at a rate of $T^{3/4}$, faster than that of $\hat{\alpha}_i$. In the model \eqref{eq:single-index general factor model} with a cointegrated single index, the convergence rates of $\hat{\alpha}_i$ are faster. Along the axis of $\alpha_{0i}/\|\alpha_{0i}\|$, the convergence rate improves to $T^{1/2}$, and along axes orthogonal to $\alpha_{0i}/\|\alpha_{0i}\|$, it improves to $T$. Consequently, $\hat{\alpha}_{i}$ achieves an overall convergence rate of $T^{1/2}$. Moreover, the convergence rate of $\hat{f}_t$ to $f_{0t}$ is $N^{1/2}$, aligning with conventional estimation rates of factors, as discussed in prior studies (c.f., \citealt{bai2003inferential}; \citealt{chen2021quantile}; \citealt{gao2023binary}). The normalized estimator $\hat{\alpha}_{i}/\|\hat{\alpha}_{i}\|$ converges to $\alpha_{0i}/\|\alpha_{0i}\|$ at an accelerated rate of $T$, surpassing the convergence rate of $\hat{\alpha}_{i}$. It is worth noting that the asymptotic distributions of the coefficient estimates are significantly different under the cointegrated single index from that under the nonstationary single index.

The modeling framework of this paper has a wide range of applications; we apply it specifically to jump arrival events in financial markets. We find that the model captures well the potential jump arrival factor, which is nonstationary. Additionally, we find that the jump arrival factor effectively explains the asset panel data, which contains information not captured by the Fama–French–Carhart five factors. The jump arrival factor---which benefits from our nonstationary binary model---is effectively applicable in empirical financial asset pricing and is distinct from jump size factors (e.g., \citealt{li2019jump}; \citealt{pelger2020understanding}).

The rest of the paper is organized as follows. In Section \ref{sec:Models and Estimation Procedure}, we introduce the model, outline the underlying assumptions, and describe the estimation procedure. Section \ref{sec:Main Results} presents the asymptotic properties of the proposed estimators. In Section \ref{sec:Simulation}, we report Monte Carlo simulation results to assess the accuracy of estimation. Section \ref{sec:Empirical Application} offers an empirical application of the model, and Section \ref{sec:Conclusion} concludes the paper.

\emph{Notation}. $\|\cdot\|$ denotes the Euclidean norm of a vector or the Frobenius norm of a matrix. For any matrix $A$ with real eigenvalues, let $\rho_{\max}(A)$ be its largest eigenvalue. Convergence in probability and in distribution are denoted by $\to_P$ and $\to_D$, respectively. $\mathbf{MN}(0,\Omega)$ denotes a mixture normal distribution with conditional covariance matrix $\Omega$. For any function $f(\cdot)$, the notation $\dot{f}(x)$, $\ddot{f}(x)$, and $\dddot{f}(x)$ refer to the first, second, and third derivatives of $f(\cdot)$ at $x$, respectively. The indicator function is written as $\mathbb{I}_{\{\cdot\}}$. We write $a_{NT}\asymp b_{NT}$ to mean that there exist positive constants $c$ and $C$, independent of $N$ and $T$, such that $c\leq a_{NT}/b_{NT}\leq C$. Similarly, $a_{NT}\lesssim b_{NT}$ means that $a_{NT}\leq Cb_{NT}$ for some positive $C$. The symbol $\otimes$ denotes the Kronecker product. The phrase ``w.p.a.1'' stands for`` with probability approaching to 1''. 

\section{Models and Estimation Procedure}\label{sec:Models and Estimation Procedure}
\subsection{Model Setup}\label{subsec:Model Setup}

Model \eqref{eq:single-index general factor model} can be rewritten in a form equivalent to 
\begin{align*}
	\left\{\begin{array}{l}
		y_{it}^* = \beta_{0i}'x_{it}+\lambda_{0i}'f_{0t}+\epsilon_{it}, \\
		y_{it} =\mathbb{I}_{\{y_{it}^*\geq0\}},
	\end{array}\right.
\end{align*}
The term  $\lambda_{0i}'f_{0t}$ captures unobserved components of individual $i$'s utility. The error term $\epsilon_{it}$ is independently and identically distributed (i.i.d.) with distribution function $\Psi:\mathbb{R}\to[0,1]$ (e.g., standard normal or standard logistic). The dependent variable $y_{it}^*$ is latent, and the observed outcome $y_{it}$ is binary, taking values of either 0 or 1. 

In our framework, both the explanatory variable $x_{it}$ and the factor $f_{0t}$ are nonstationary processes, integrated of order one (denoted as $I(1)$). We first assume that there exist neighborhoods around the true parameters $\beta_{0i}$ and $\lambda_{0i}$ such that $\beta_i$ and $\lambda_i$ always lie within these neighborhoods, ensuring that $\beta_{i}'x_{it}$ and $\lambda_{i}'f_{0t}$ remain $I(1)$ processes. In Section \ref{subsec:Theoretical Results for Nonstationary Models}, we treat $\beta_{0i}'x_{it}+\lambda_{0i}'f_{0t}$ as an $I(1)$ process. In Section \ref{subsec:Cointegrated Single-Indexes}, we allow for cointegration between $\beta_{0i}'x_{it}$ and $\lambda_{0i}'f_{0t}$, i.e., $\beta_{0i}'x_{it}+\lambda_{0i}'f_{0t}\sim I(0)$. These two settings encompass both stationary and nonstationary single indices and are broadly applicable.

The assumptions regarding $x_{it}$ and $f_{0t}$ are outlined below for the subsequent development of the theory.
\begin{assumption}\label{assump:I(1) process}(Integrated Processes)
\begin{enumerate}[(i)]
	\item The covariate series follow $x_{it}=x_{it-1}+e_{it}$ and the factor series follow $f_{0t}=f_{0t-1}+v_{t}$, with initial conditions $x_{i0}=O_P(1)$ and $f_{00}=O_P(1)$, where $\{e_{it}\}$ and $\{v_t\}$ are linear processes defined as $e_{it}=\varPi^e(L)\varepsilon_{it}^{e}=\sum_{k=0}^{\infty}\varPi_k^e\varepsilon_{it-k}^{e}$ and $v_{t}=\varPi^v(L)\varepsilon_{t}^{v}=\sum_{k=0}^{\infty}\varPi_k^v\varepsilon_{t-k}^{v}$, where $\varPi^j(L)=\sum_{k=0}^{\infty}\varPi_k^j L^k$ with the lag operator $L$ for $j\in\{e,v\}$, $\varPi^e(1)$ and $\varPi^v(1)$ are nonsingular, and $\varPi_0^e=I_{q}$ and $\varPi_0^v=I_{r}$. Additionally, we assume $\sum_{k=0}^{\infty}k(\|\varPi_k^e\|+\|\varPi_k^v\|)<\infty$. The vector $\{\varepsilon_{it}=(\varepsilon_{it}^{e'},\varepsilon_{t}^{v'})'\}$ are i.i.d. with mean zero and satisfies $E\|\varepsilon_{it}\|^\iota<\infty$ for some $\iota>8$. The distribution of $(\varepsilon_{it})$ is absolutely continuous with respect to the Lebesgue measure and has a characteristic function $\varphi_i(t)$ such that $\varphi_i(t)=o(\|t\|^{-\kappa})$ as $\|t\|\to\infty$ for some $\kappa>0$.
	\item Define $U_{i,n}(s)=\frac{1}{\sqrt{T}}\sum_{t=1}^{[Ts]}u_{it}$, $E_{i,n}(s)=\frac{1}{\sqrt{T}}\sum_{t=1}^{[Ts]}e_{it}$, and $V_{n}(s)=\frac{1}{\sqrt{T}}\sum_{t=1}^{[Ts]}v_{t}$. We assume that there exist $(q+r+1)$-dimensional Brownian motion $(U_i,E_i,V)$ such that $(U_{i,n}(s),E_{i,n}(s),V_{n}(s))\to_D(U_i(s),E_i(s),V(s))$ in $D[0,1]^{q+r+1}$ for all $i\in\{1,...,N\}$.
\end{enumerate}
\end{assumption}

These conditions are standard in nonstationary model estimation. In particular, Assumption \ref{assump:I(1) process}(ii) is widely used in related studies, such as \cite{park2000nonstationary}, \cite{dong2016estimation}, and \cite{trapani2021inferential}. We define $H_i(s)=(E_i(s)',V(s)')'$.

\begin{assumption}\label{assump:Factors and loadings}(Covariate Coefficients, Factors, and Factor Loadings)
	\begin{enumerate}[(i)]
		\item There exists a positive constant $C$, such that $\|\beta_{0i}\|\leq C$ and $\|\lambda_{0i}\|\leq C$ for all $i=1,...,N$.
		\item As $T\to\infty$, $\sum_{t=1}^Tf_{0t}f_{0t}'/T^2\to_D\int_0^1W_sW_s'ds$, where $W_s$ is a vector of Brownian motions with a positive definite covariance matrix.
		\item $\sum_{i=1}^N\lambda_{0i}\lambda_{0i}'/N=\text{diag}(\sigma_{N1},...,\sigma_{Nr})$ with $\sigma_{N1}\geq\sigma_{N2}\geq\sigma_{Nr}$, and $\sigma_{Ni}\to\sigma_i$ as $N\to\infty$ for $i=1,...,r$, where $\infty>\sigma_1>\sigma_2>\cdots>\sigma_r>0$.
	\end{enumerate}
\end{assumption}
Assumption \ref{assump:Factors and loadings}(i) means that the loadings are in compact set. Sufficient conditions for Assumption \ref{assump:Factors and loadings}(ii) can be found in \cite{hansen1992convergence}. The assumption of positive definiteness in Assumption \ref{assump:Factors and loadings}(ii) precludes cointegration among the components of $f_{0t}$. Similar assumptions are made in \cite{bai2004estimating}. Assumption \ref{assump:Factors and loadings}(iii) is a version of the strong factors assumption, which is commonly used in the literature, such as in \cite{bai2003inferential} and \cite{chen2021quantile}. The requirement that $\sigma_1,...,\sigma_r$ are distinct is similar to the assumption in \cite{bai2003inferential}, which provides a convenient way to identify and order the factors.

\subsection{Estimation Procedure}\label{subsec:Estimation Procedure}

Let $z_{it}=\beta_{i}'x_{it}+\lambda_{i}'f_{t}$, $\alpha_{i}=(\beta_{i}',\lambda_{i}')'$, $g_{it}=(x_{it}',f_{t}')'$, $A=(\alpha_{1},...,\alpha_{N})'$, $B=(\beta_{1},...,\beta_{N})'$, $\Lambda=(\lambda_{1},...,\lambda_{N})'$, and $F=(f_{1},...,f_{T})'$. Let $g_{0it}$, $A_0$, $B_0$, $\Lambda_0$, and $F_0$ be the true parameters. 

Given the stochastic properties of $\{u_{it}\}$, the log-likelihood function is
\begin{align}\label{eq:log-likelihood function}
	\log L(B,\Lambda,F)=\sum_{i=1}^N\sum_{t=1}^T\left[y_{it}\log\Psi(z_{it})+(1-y_{it})\log(1-\Psi(z_{it}))\right].
\end{align}
Due to the rotational indeterminacy of the factor loadings $\lambda_i$ and the factors $f_t$, we impose identification constraints, following standard approaches:
\begin{align}\label{eq:identification conditions}
		\begin{aligned}
	&\mathcal{F}=\left\{F\in\mathbb{R}^{T\times r}:\frac{F'F}{T^2}=I_r\right\},\\
	&\mathcal{L}=\left\{\Lambda\in\mathbb{R}^{N\times r}:\frac{\Lambda'\Lambda}{N}\text{is diagonal with non-increasing elements}\right\}.
	\end{aligned}
\end{align}

The estimators for $(B_0,\Lambda_0,F_0)$ are then defined as:
\begin{align}\label{eq:MLE estimators}
	(\hat{B},\hat{\Lambda},\hat{F})=\argmin_{B\in\mathbb{R}^{N\times q},\Lambda\in\mathcal{L},F\in\mathcal{F}}\log L(B,\Lambda,F).
\end{align}
Direct optimization of this expression is challenging due to the absence of closed-form solutions for our estimators, unlike in principal component analysis (PCA). To address this, we propose an iterative optimization algorithm.

Define the serial and cross-section averages as
\begin{align*}	 \mathbb{L}_{iT}(\alpha_i,F)=&\sum_{t=1}^Tl_{it}(z_{it}),\quad\text{and}\quad \mathbb{L}_{Nt}(A,f_t)=\sum_{i=1}^Nl_{it}(z_{it}),
\end{align*}
where $l_{it}(z_{it})=\left[y_{it}\log\Psi(z_{it})+(1-y_{it})\log(1-\Psi(z_{it}))\right]$. The iterative optimization algorithm proceeds as follows:
\textit{\begin{enumerate}[~~~~Step 1:]
	\item Randomly select the initial parameter $A^{(0)}$.
	\item Given $A^{(l-1)}$, solve $f_t^{(l-1)}=\argmin_{f}\mathbb{L}_{Nt}\left(A^{(l-1)},f\right)$ for $t=1,...,T$; given $F^{(l-1)}$, solve $\alpha_i^{(l)}=\argmin_{\alpha}\mathbb{L}_{iT}\left(\alpha,F^{(l-1)}\right)$ for $i=1,...,N$.
	\item Repeat Step 2 until a convergence criterion is met.
	\item Let $\Lambda^*$ and $F^*$ be the final estimators after the iteration process. Normalize $\Lambda^*$ and $F^*$ to satisfy the constraints in \eqref{eq:identification conditions}.
\end{enumerate}}

In Step 3, various tolerance conditions can be employed, such as those based on parameter changes or the objective function. In this paper, we terminate the iteration when the change in the objective function is below a threshold $\varrho$, i.e., $|L^{\text{new}}-L^{\text{last time}}|<\varrho$, where $L^{\text{new}}$ and $L^{\text{last time}}$ denote the current and previous values of the  log-likelihood function, respectively. In Step 4, we compute
\begin{align*}
	 \left(\frac{F^{*'}F^*}{T^2}\right)^{1/2}\left(\frac{\Lambda^{*'}\Lambda^*}{N}\right)\left(\frac{F^{*'}F^*}{T^2}\right)^{1/2}=Q^*D^*Q^*,
\end{align*}
where $D^*$ is a diagonal matrix. We then sort the diagonal elements of $\Lambda^*(F^{*'}F^*/T^2)^{1/2}Q^*$ in descending order to obtain $\hat{\Lambda}$. Similarly, we sort $F^*(F^{*'}F^*/T^2)^{-1/2}Q^*$ according to the same order to obtain $\hat{F}$.

\begin{remark}
	If the covariate component $\beta_{0i}x_{it}$ is absent, the binary probability simplifies to $\Psi(\lambda_{0i}'f_{0t})$, reducing the model to a pure binary factor model. Although the likelihood function remains non-convex in $(\Lambda,F)$, the limit of $\mathbb{L}_{iT}(\lambda_i,F)$ becomes globally convex for each $\lambda_i$ given $F$, and the limit of $\mathbb{L}_{Nt}(\Lambda,f_t)$ becomes globally convex for each $f_t$ given $\Lambda$. Consequently, these two optimization problems can be efficiently solved, as demonstrated in prior research, e.g., \cite{chen2021quantile}. To illustrate this, consider the logit model. The leading terms in the Hessian matrices are:
	\begin{align*}
		\frac{\partial \mathbb{L}_{iT}(\lambda,F)}{\partial\lambda\partial\lambda'}=-\sum_{t=1}^T\eta_{it}(1-\eta_{it})f_{0t}f_{0t}',\quad\frac{\partial \mathbb{L}_{Nt}(\Lambda,f)}{\partial f\partial f'}=-\sum_{i=1}^N\eta_{it}(1-\eta_{it})\lambda_{0i}\lambda_{0i}',
	\end{align*}
where $\eta_{it}(1-\eta_{it})=\frac{e^{\lambda_{0i}'f_{0t}}}{(1+e^{\lambda_{0i}'f_{0t}})^2}$ is the logistic function, a nonlinear integrable function of $f_{0t}$.

When appropriately normalized, the Hessian matrix $\frac{\partial \mathbb{L}_{iT}(\lambda,F)}{\partial\lambda\partial\lambda'}$ converges weakly to a random limit matrix, rather than a constant matrix (see \citealt{park2000nonstationary}). Additionally, the Hessian matrix $\frac{\partial \mathbb{L}_{Nt}(\Lambda,f)}{\partial f\partial f'}$
may converge to a neighborhood of zero when $t$ is large (due to $x\mapsto e^{x}/(1+e^{x})^2$ is integrable and $\frac{e^{\lambda_{0i}'f_{0t}}}{(1+e^{\lambda_{0i}'f_{0t}})^2}\approx0$ outside the effective range of this function), making the consistency of $f_t$ more challenging to ensure. However, when $\frac{\partial \mathbb{L}_{Nt}(\Lambda,f)}{\partial f\partial f'}$ is adjusted with respect to $t$ (see Section \ref{sec:Main Results}), it also converges weakly to a stochastic limit matrix. Both Hessian matrices are almost surely negative definite, ensuring that the limit functions are globally concave.
\end{remark}

\section{Main Results}\label{sec:Main Results}

In this section, we delve into the theoretical foundations of the maximum likelihood estimation presented in \eqref{eq:MLE estimators}. To facilitate understanding, we introduce the class of regular functions.
\begin{definition}\label{def:function class}
A function $f\in\mathbb{F}_R:\mathbb{R}\to\mathbb{R}$ is termed \emph{regular} if it satisfies the following conditions: (i) $|f(x)|\leq M$ for some $M>0$ and all $x\in\mathbb{R}$; (ii) $\int_{-\infty}^{\infty}|f(x)|dx<\infty$; (iii) $f$ is differentiable with a bounded derivative. Let $\mathbb{F}_I$ be the set of functions that are bounded and integrable, and $\mathbb{F}_B$ be the set of bounded functions that vanish at infinity. Clearly, the inclusions $\mathbb{F}_R\subset\mathbb{F}_I\subset\mathbb{F}_B$ hold.
\end{definition}

Next, we define the leading terms in score and Hessian as follows.
\begin{align*}
	M=\frac{\dot{\Psi}}{\Psi(1-\Psi)},\quad K=M\dot{\Psi}=M^2\Psi(1-\Psi).
	\end{align*}
For the Logit case, $M(x)=1$ and $K(x)=e^x/(1+e^x)^2$. For the Probit case, $M(x)=\phi(x)/(\Phi(x)(1-\Phi(x)))$ and $K(x)=\phi^2(x)/(\Phi(x)(1-\Phi(x)))$, where $\phi$ is the probability density function of standard normal variable, and $\Phi$ is the cumulative distribution function of the standard normal distribution.

\subsection{Theoretical Results for Nonstationary Models}\label{subsec:Theoretical Results for Nonstationary Models}

In this subsection, we assume that the process $\{g_{0it}\}$ is integrated and of full rank, indicating the absence of cointegrating relationships among its component time series. To develop an asymptotic theory for the estimators $(\hat{A},\hat{F})$ defined in \eqref{eq:MLE estimators}, we impose several assumptions on the functions $\Psi$, $M$, and $K$.
\begin{assumption}\label{assump:function categories}(Function Categories)
The function $\Psi$ is three times differentiable on $\mathbb{R}$. Additionally, the following conditions hold: (i) $K_2\in\mathbb{F}_R$; (ii) $\dot{\Psi}_1$, $(\dot{M}\dot{\Psi})$, $(M\ddot{\Psi})_2$, $\ddot{M}\dot{\Psi}$, $(\ddot{M}\Psi^{1/2}(1-\Psi)^{1/2})_2\in\mathbb{F}_I$; (iii) $(\dot{M}\Psi^{1/2}(1-\Psi)^{1/2})_2$, $(M^3\dot{\Psi})_4\in\mathbb{F}_B$.
\end{assumption}
These assumptions are mild and are satisfied by common models such as logit and probit. For notational convenience in the subsequent derivations, we define $F_v=(f_1',...,f_T')'$ and $A_v=(\alpha_1',...,\alpha_N')'$. 

To facilitate the analysis, we introduce block matrices for the score and Hessian. Define $S_{NT,1}(A,F)=\frac{\partial}{\partial A_v}\log L(A,F)$, $J_{NT,11}(A,F)=\frac{\partial^2}{\partial A_v\partial A_v'}\log L(A,F)$, $S_{NT,2}(A,F)=\frac{\partial}{\partial F_v}$ $\log L(A,F)$, and $J_{NT,22}(A,F)=\frac{\partial^2}{\partial F_v\partial F_v'}\log L(A,F)$. The score function with respect to $A$ is expressed as $S_{NT,1}(A,F)=\left(\left(S_{NT,1}^{(1)'}(\alpha_1,F),...,S_{NT,1}^{(N)'}(\alpha_N,F)\right)'\right)_{N(q+r)\times 1}$. The corresponding Hessian is $J_{NT,11}(A,F)=\left(\mathrm{diag}\left(J_{NT,11}^{(1)}(\alpha_1,F),...,J_{NT,11}^{(N)}(\alpha_N,F)\right)\right)_{N(q+r)\times N(q+r)}$. Similarly, the score function with respect to $F$ is $S_{NT,2}(A,F)=\left(\left(S_{NT,2}^{(1)'}(A,f_1),...,S_{NT,2}^{(T)'}(A,f_T)\right)'\right)_{Tr\times 1}$, and the Hessian is $J_{NT,22}(A,F)=\left(\mathrm{diag}\left(J_{NT,22}^{(1)}(A,f_1),...,J_{NT,22}^{(T)}(A,f_T)\right)\right)_{Tr\times Tr}$. The diagonal structure of the Hessian is straightforward to verify. For $\alpha_i$,
\begin{align*} S_{NT,1}^{(i)}(\alpha_i,F)=&\sum_{t=1}^TM(z_{it})g_{it}(y_{it}-\Psi(z_{it})),\\ J_{NT,11}^{(i)}(\alpha_i,F)=&-\sum_{t=1}^TK(z_{it})g_{it}g_{it}'+\sum_{t=1}^T\dot{M}(z_{it})g_{it}g_{it}'(y_{it}-\Psi(z_{it})).
\end{align*}
For $f_t$, 
\begin{align*} S_{NT,2}^{(t)}(A,f_t)=&\sum_{i=1}^NM(z_{it})\lambda_{i}(y_{it}-\Psi(z_{it})),\\ J_{NT,22}^{(t)}(A,f_t)=&-\sum_{i=1}^NK(z_{it})\lambda_{i}\lambda_{i}'+\sum_{i=1}^N\dot{M}(z_{it})\lambda_{i}\lambda_{i}'(y_{it}-\Psi(z_{it})).
\end{align*}

Due to the unit root behavior of ${z_{0it}}$, its probability mass spreads out in a manner similar to a Lebesgue type. Given that $K$ is integrable (as per Assumption \ref{assump:function categories}), $K(z_{0it})\approx0$ outside the effective range of $K$. This indicates that only moderate values of $z_{0it}$ prevent $K(z_{0it})$ from diminishing. Unlike $\sum_{t=1}^TK(z_{0it})g_{0it}g_{0it}'$, the sum $\sum_{i=1}^NK(z_{0it})\lambda_{0i}\lambda_{0i}'$ varies with time $t$, reflecting the spread of $z_{0it}$ at specific time points. Therefore, to normalize the Hessian appropriately, it is crucial to analyze the convergence behavior of $\sum_{i=1}^N K(z_{0it})$ and select a suitable normalizing sequence concerning $t$. Based on this analysis, we introduce the following assumptions.

\begin{assumption}\label{assump:some integrable functions}(Integrable Functions)
		\begin{enumerate}[(i)]
		\item For a normally distributed random variable $z\sim\mathbf{N}(0,\sigma^2)$, we assume $E(K(z))\asymp\sigma^{-2\delta}$. Additionally,  $E(\dot{K}(z)z)\vee E(\ddot{M}(\Psi(1-\Psi))^{1/2}(z))\vee E(K^2(z)z^2)\vee E(K^2(z))\lesssim\sigma^{-2\delta}$ for $\sigma^{-1}=o(1)$ and $\delta\in(1/4,3/4)$.
		\item The variance $Var\left(\frac{t^{\delta}}{N}\sum_{i=1}^NK(z_{0it})\right)\to0$ as $N\to \infty$ for all $t=1,...,T$.
		\item For each $t$, as $N\to\infty$, $\frac{t^{\delta/2}}{\sqrt{N}}\sum_{i=1}^NM(z_{0it})\lambda_{0i}u_{it}\to_D\mathbf{N}(0,\Omega_{f,t})$, where the covariance matrix $\Omega_{f,t}=\mathrm{lim}_{N\to\infty}\frac{t^{\delta}}{N}\sum_{i=1}^NE\left[K(z_{0it})\right]\lambda_{0i}\lambda_{0i}'$.
		\item The sequences satisfy $T^{\delta'}/N=o(1)$ and $N^{\delta''}/T=o(1)$ for some $\delta''>0$, where $\delta'\geq\max(\delta,\delta/2+1/2)$.
		\item For functions $f$ defined in Assumptions \ref{assump:function categories}(i) and (ii),  $\frac{1}{\sqrt{T}}\sum_{t=1}^Tt^{\delta}f\left(h_{0it}^{(1)}\right)=O_P(1)$, where $h_{0it}^{(1)}=z_{0it}/\|a_{0i}\|$.
	\end{enumerate}
\end{assumption}
Assumption \ref{assump:some integrable functions}(i) is mild. To illustrate its plausibility, simulations (reported in the Supplementary Material) indicate that $\delta$ is approximately 0.5 in both the logit and probit models. Assumption \ref{assump:some integrable functions}(ii) ensures that the sum $\sum_{t=1}^TK(z_{0it})$ appropriately utilizes the properties outlined in Assumption \ref{assump:some integrable functions}(i). In the degenerate case, where $z_{0it}$ is independent across $i$ and $t$ takes moderate values, the variance satisfies $Var\left(\frac{t^{\delta}}{N}\sum_{i=1}^NK(z_{0it})\right)=O(1/N)$. In the stationary case, the $\delta$ in Assumption \ref{assump:some integrable functions}(iii) should be 0. Assumption \ref{assump:some integrable functions}(iv) imposes a constraint on the sample size and the dimensionality. Assumption \ref{assump:some integrable functions}(v) addresses information overflow in the cross-section and therefore requires a slightly stronger condition than Assumption \ref{assump:function categories}.

As is standard, we have the following Taylor expansions.
\begin{align}\label{eq:Taylor expansion for hat}
	\begin{aligned}
	 0=S_{NT,1}^{(i)}(\hat{\alpha}_i,\hat{F})=&S_{NT,1}^{(i)}(\alpha_{0i},F_0)+J_{NT,11}^{(i)}(\alpha_i,F)(\hat{\alpha}_i-\alpha_{0i})+\sum_{t=1}^TJ_{NT,12}^{(i,t)}(\alpha_i,f_t)(\hat{f}_t-f_{0t}),\\
	 0=S_{NT,2}^{(t)}(\hat{A},\hat{f}_t)=&S_{NT,2}^{(t)}(A_{0i},f_{0t})+J_{NT,22}^{(t)}(A_{i},f_{t})(\hat{f}_t-f_{0i})+\sum_{i=1}^NJ_{NT,21}^{(t,i)}(\alpha_i,f_t)(\hat{\alpha}_i-\alpha_{0t}),
	\end{aligned}
	\end{align}
where $J_{NT,12}^{(i,t)}$ denotes the $(i,t)$th block of the matrix $J_{NT,12}=\frac{\partial^2}{\partial A_v\partial F_v'}\log L(A,F)$, and $J_{NT,21}^{(t,i)}$ denotes the $(t,i)$th block of the matrix $J_{NT,21}=\frac{\partial^2}{\partial F_v\partial A_v'}\log L(A,F)$. Additionally, $(\alpha_i',f_t')$ is some point between $(\hat{\alpha}_i',\hat{f}_t')$ and $(\alpha_{0i}',f_{0t}')$.

The asymptotic theory for $(\hat{A},\hat{F})$ can be derived from Equation \eqref{eq:Taylor expansion for hat}. To aid in the development of this theory, we rotate the coordinate system based on the true parameter $A_0$ using an orthogonal matrix $Q_i\in\mathbb{R}^{(q+r)\times(q+r)}$, where $Q_i=(Q_i^{(1)},Q_i^{(2)})$ and $Q_i^{(1)}=\alpha_{0i}/\|\alpha_{0i}\|$.\footnote{Because estimators converge at different rates in both the parallel and orthogonal directions relative to $\alpha_{0i}$, performing a rotation enables a more comprehensive theoretical examination.} This matrix $Q_i$ will be used to rotate all vectors in $\mathbb{R}^{q+r}$ for $i=1,...,N$. Specifically, we define the following quantities:
\begin{align*}
	 &\theta_{0i}:=Q_i'\alpha_{0i}=(\theta_{0i}^{(1)},\theta_{0i}^{(2)'})'\quad\text{where}\quad\theta_{0i}^{(1)}=\|\alpha_{0i}\|,\theta_{0i}^{(2)}=Q_{i}^{(2)'}\alpha_{0i}=0,\\
	 &h_{0it}:=Q_i'g_{0it}=(h_{0it}^{(1)},h_{0it}^{(2)'})'\quad\text{where}\quad h_{0it}^{(1)}=\alpha_{0i}'g_{0it}/\|\alpha_{0i}\|=z_{0it}/\|\alpha_{0i}\|,h_{0it}^{(2)}=Q_{i}^{(2)'}g_{0it}.
\end{align*}
In the general case, we define $\theta_{i}:=Q_i'\alpha_{i}$ and $h_{it}:=Q_i'g_{it}$. With these definitions, we can rewrite the model as $y_{it}=\Psi(\alpha_{0i}'Q_iQ_i'g_{0it})+u_{it}=\Psi(\theta_{0i}'h_{0it})+u_{it}$. By Assumption \ref{assump:I(1) process} and applying the continuous mapping theorem, we obtain the following convergence results for $s\in[0,1]$,
\begin{align*}
	 \frac{1}{\sqrt{T}}h_{0i[Ts]}^{(1)}\to_DH_{1i}(s)=Q_i^{(1)'}H_i(s)\quad\text{and}\quad\frac{1}{\sqrt{T}}h_{0i[Ts]}^{(2)}\to_DH_{2i}(s)=Q_i^{(2)'}H_i(s).
\end{align*}
It is important to note that the rotation is not required in practice, and indeed, is conceptually impossible since $A_0$ is unknown. The rotation serves only as a tool for deriving the asymptotic theory for the proposed estimators. 

If $\hat{\theta}_i$ is the maximum likelihood estimator of $\theta_{0i}$, then we have $\hat{\theta}_i=Q_i'\hat{\alpha}_i$. The score function and Hessian for the parameter $\theta_i$ can be expressed in terms of $\alpha_i$ as follows: $S_{NT,1}^{(i)}(\theta_i,F)=Q_i'S_{NT,1}^{(i)}(\alpha_i,F)$ and $J_{NT,11}^{(i)}(\theta_i,F)=Q_i'J_{NT,11}^{(i)}(\alpha_i,F)Q_i'$. Using this relationship, we can derive the following Taylor expansion:
\begin{align}\label{eq:Taylor expansion for hat of re-parameter}
	\begin{aligned}
		 0=S_{NT,1}^{(i)}(\hat{\theta}_i,\hat{F})=&S_{NT,1}^{(i)}(\theta_{0i},F_0)+J_{NT,11}^{(i)}(\theta_i,F)(\hat{\theta}_i-\theta_{0i})+\sum_{t=1}^TJ_{NT,12}^{(i,t)}(\theta_i,f_t)(\hat{f}_t-f_{0t}).
	\end{aligned}
\end{align}
Define $\Theta_0=(\theta_{01},...,\theta_{0N})'$, $D_T=\mathrm{diag}(T^{1/4},T^{3/4}I_{q+r-1})$, $B_N=\mathrm{diag}(1^{\delta/2},2^{\delta/2},...,T^{\delta/2})$ and $C_{NT}=\min\{\sqrt{N},\sqrt{T}\}$.

\begin{assumption}\label{assump:max eigenvalues for covariance}(Covariance)
	$\rho_{\max}(\mathcal{A}_{11})$, $\rho_{\max}(\mathcal{A}_{22}),$ $\rho_{\max}(\mathcal{A}_{11}^{-1})$, $\rho_{\max}(\mathcal{A}_{22}^{-1})$, $\rho_{\max}(\mathcal{A}_{12}\mathcal{A}_{22}^{-1}\mathcal{A}_{21})$, and $\rho_{\max}(\mathcal{A}_{21}\mathcal{A}_{11}^{-1}\mathcal{A}_{12})$ are finite, where 
	$\mathcal{A}_{11}=(I\otimes D_T)^{-1}\mathcal{J}_{NT,11}(1)(I\otimes D_T)^{-1}$, $\mathcal{A}_{12}=(I\otimes D_T)^{-1}\mathcal{J}_{NT,12}(1)\frac{(B_N\otimes I)}{\sqrt{N}}$, $\mathcal{A}_{21}=\mathcal{A}_{21}'$, $\mathcal{A}_{22}=\frac{(B_N\otimes I)}{\sqrt{N}}\mathcal{J}_{NT,22}(1)\frac{(B_N\otimes I)}{\sqrt{N}}$, with the $i$th block of $\mathcal{J}_{NT,11}(1)$ given by $\mathcal{J}_{NT,11}^{(i)}(1)=-\sum_{t=1}^TK(z_{0it})h_{0it}h_{0it}'$, the $t$th block of $\mathcal{J}_{NT,22}(1)$ given by $\mathcal{J}_{NT,22}^{(t)}(1)=-\sum_{i=1}^NK(z_{0it})\lambda_{0i}\lambda_{0i}'$, and the $(i,t)$th block of $\mathcal{J}_{NT,12}(1)$ given by $\mathcal{J}_{NT,12}^{(i)}(1)=-K(z_{0it})h_{0it}\lambda_{0i}$.
\end{assumption}
Assumption \ref{assump:max eigenvalues for covariance} is used to deriving results related to the inverse of the Hessian matrix. It is a mild assumption because both $\mathcal{A}_{11}$ and $\mathcal{A}_{22}$ are diagonal matrices.

We now present the average rate of convergence for $\hat{\Theta}$ and $\hat{F}$.
\begin{theorem}\label{thm:average convergence rate}
	 Under Assumptions \ref{assump:I(1) process}-\ref{assump:max eigenvalues for covariance}, the following results hold:
	\begin{align*}
		 \frac{1}{\sqrt{N}}\|T^{-1/4}(\hat{\Theta}-\Theta_0)D_T\|=O_P\left(T^{1/4}C_{NT}^{-1}\right)\quad\text{and}\quad\frac{1}{\sqrt{T}}\|B_N^{-1}(\hat{F}-F_0)\|=O_P\left(C_{NT}^{-1}\right).
	\end{align*}
\end{theorem}
Theorem \ref{thm:average convergence rate} shows that the estimator $\hat{\Theta}$ exhibits dual convergence rates.
\begin{enumerate}[(i)]
	\item Along the coordinates parallel to $\{\alpha_{0i},i=1,...,N\}$ (i.e., $\left(\hat{\theta}_{1}^{(1)},...,\hat{\theta}_{N}^{(1)}\right)'/\sqrt{N}$), the average rate of convergence is $T^{1/4}C^{-1}_{NT}$.
	\item Along the coordinates orthogonal to $\{\alpha_{0i},i=1,...,N\}$ (i.e., $\left(\hat{\theta}_{1}^{(2)},...,\hat{\theta}_{N}^{(2)}\right)'/\sqrt{N}$), the average rate of convergence is $T^{-1/4}C_{NT}^{-1}$.
\end{enumerate}
For $\hat{F}$, the rate of convergence depends on $B_N$, with the difference $\hat{f}_t-f_{0t}$ scaled by $t^{\delta/2}$. If $T^{\delta}/N=o(1)$, the estimators for $\{\hat{f}_t\}$ are consistent. Otherwise, only a subset of $\{\hat{f}_t\}$ are consistent. It is noteworthy that the estimation of $f_{0t}$ becomes less accurate for larger $t$. This is intuitive, as $Var(f_{0t})=O(t)$, meaning that the uncertainty increases as $t$ grows. To the best of our knowledge, this is the first time such a phenomenon has been observed in the context of factor estimators. The explanation lies in the fact that $J_{NT,22}^{(t)}(A_0,f_{0t})\approx0$ for larger $t$.

After performing a rotation $A=(Q_1\theta_1,...,Q_N\theta_N)$, we present the convergence rates for $\hat{A}$ in Corollary 3.2 below.
\begin{corollary}\label{coro:average convergence rate}
	Under Assumptions \ref{assump:I(1) process}-\ref{assump:max eigenvalues for covariance},
		\begin{align*}
		 \frac{1}{\sqrt{N}}\|(\hat{A}-A_0)\|=O_P\left(T^{1/4}C_{NT}^{-1}\right).
	\end{align*}
\end{corollary}
Corollary \ref{coro:average convergence rate} demonstrates the collective consistency of the estimator $\hat{A}$. The subsequent theorem provides the asymptotic distributions of $\hat{\theta}_i$ and $\hat{f}_t$.

\begin{theorem}\label{thm:asymptotic distribution for theta and f}
	Under Assumptions \ref{assump:I(1) process}-\ref{assump:max eigenvalues for covariance}, as $N,T\to\infty$,
	\begin{align*}
		 D_T(\hat{\theta}_i-\theta_{0i})\to_D\Omega_{\theta,i}^{-1/2}W_i(1)\quad\text{and}\quad \sqrt{N}t^{-\delta/2}(\hat{f}_t-f_{0t})\to_D\mathbf{N}(0,\Omega_{f,t}^{-1}),
	\end{align*}
where $W_i(1)$ is $(q+r)$-dimensional vector of Brownian motion with covariance matrix $I_{q+r}$ and independent of $H$. Let
\begin{align*}
	\Omega_{\theta,i}=
	\begin{pmatrix}
		L_{1i}(1,0)\int_{\mathbb{R}}m^2K(\|\alpha_{0i}\|m)dm & \int_{0}^1H_{2i}(s)'dL_{1i}(s,0)\int_{\mathbb{R}}mK(\|\alpha_{0i}\|m)dm \\
		\int_{0}^1H_{2i}(s)dL_{1i}(s,0)\int_{\mathbb{R}}mK(\|\alpha_{0i}\|m)dm & \int_{0}^1H_{2i}(s)H_{2i}(s)'dL_{1i}(s,0)\int_{\mathbb{R}}K(\|\alpha_{0i}\|m)dm
	\end{pmatrix}
\end{align*}
and $\Omega_{f,t}=\mathrm{lim}_{N\to\infty}\frac{t^{\delta}}{N}\sum_{i=1}^NE\left[K(z_{0it})\right]\lambda_{0i}\lambda_{0i}'$, where $L_{1i}(s,0)=L_{H_{1i}}(s,0)\sigma_{H_{1i}}$ with $L_{H_{1i}}(s,0)$ being the local time of $H_{1i}$ and $\sigma_{H_{1i}}$ its variance.
\end{theorem}

The asymptotic behavior of the estimator $\hat{f}_t$ varies with $t$, influenced by the Hessian matrix. Recall the notation of $D_T$, two distinct limiting distributions for $\hat{\theta}_i$ emerge.
\begin{align}\label{eq:asymptotic limits of theta1 and theta2}
	 T^{1/4}\left(\hat{\theta}_{i}^{(1)}-\theta_{0i}^{(1)}\right)\to_D\mathbf{MN}\left(0,\bar{\omega}_{11}^{\theta,i}\right)\quad\text{and}\quad T^{3/4}\left(\hat{\theta}_{i}^{(2)}-\theta_{0i}^{(2)}\right)\to_D\mathbf{MN}\left(0,\bar{\omega}_{22}^{\theta,i}\right),
\end{align}
where 
\begin{align*} \bar{\omega}_{11}^{\theta,i}=\left(\omega_{11}^{\theta,i}-\omega_{12}^{\theta,i}(\omega_{22}^{\theta,i})^{-1}\omega_{21}^{\theta,i}\right)^{-1}\quad\text{and}\quad\bar{\omega}_{22}^{\theta,i}=\left(\omega_{22}^{\theta,i}-\omega_{21}^{\theta,i}(\omega_{11}^{\theta,i})^{-1}\omega_{12}^{\theta,i}\right)^{-1},
\end{align*}
with $\Omega_{\theta,i}:=\begin{pmatrix} \omega_{11}^{\theta,i}& \omega_{12}^{\theta,i}\\\omega_{21}^{\theta,i}&\omega_{22}^{\theta,i} \end{pmatrix}$ and $\omega_{11}^{\theta,i}=L_{1i}(1,0)\int_{\mathbb{R}}m^2K(\|\alpha_{0i}\|m)dm$. 

The dual convergence rates presented in Equation \eqref{eq:asymptotic limits of theta1 and theta2} are not surprising; similar results have been observed in various problems involving nonlinear functions, such as \cite{park2000nonstationary} and \cite{dong2016estimation}. This implies that, in multivariate cases ($q+r>1$), modest values of $\{g_{0it}\}$ significantly influence a nonlinear function along $\alpha_{0i}/\|\alpha_{0i}\|$. In contrast, there are no such restrictions on $\{g_{0it}\}$ in the direction orthogonal to $\alpha_{0i}/\|\alpha_{0i}\|$, allowing larger values of $\{g_{0it}\}$ to contribute.

We introduce the normalized estimators $\hat{\alpha}_{i}^{\circ}:=\hat{\alpha}_i/\|\hat{\alpha}_i\|$ and  $\hat{\theta}_{i}^{\circ}:=\hat{\theta}_i/\|\hat{\theta}_i\|$, derived from $\hat{\theta}_i=Q_i\hat{\alpha}_i$ and $\|\hat{\alpha}_i\|=\|\hat{\theta}_i\|$. Specifically, $\hat{\theta}_{i}^{\circ}=(\hat{\theta}_{i}^{(1)\circ},\hat{\theta}_{i}^{(2)\circ'})':=(\hat{\theta}_{i}^{(1)}/\|\hat{\theta}_i\|,\hat{\theta}_{i}^{(2)'}/\|\hat{\theta}_i\|)'$. The following corollary characterizes the asymptotic behavior of $\hat{\theta}_{i}^{(1)\circ}$ and $\hat{\theta}_{i}^{(2)\circ}$.
\begin{corollary}\label{coro:normalized CLT}
	Under Assumptions \ref{assump:I(1) process}-\ref{assump:max eigenvalues for covariance}, as $T\to\infty$
	\begin{align*}
		 T^{3/2}(\hat{\theta}_{i}^{(1)\circ}-1)\to_D-\frac{1}{2\|\alpha_{0i}\|^2}\|\zeta_i\|^2\quad\text{and}\quad T^{3/4}\hat{\theta}_{i}^{(2)\circ}\to_D\frac{1}{\|\alpha_{0i}\|}\zeta_i,
	\end{align*}
where $\zeta_i\sim\mathbf{MN}(0,\bar{\omega}_{22}^{\theta,i})$.
\end{corollary}

After normalization, the convergence rate along the $\alpha_{0i}/\|\alpha_{0i}\|$ direction increases to $T^{3/2}$, while in the orthogonal direction, it remains at $T^{3/4}$. By leveraging the linear relationship between $\hat{\alpha}_i$ and $\hat{\theta}_i$ (and similarly between $\hat{\alpha}_i^{\circ}$ and $\hat{\theta}_i^{\circ}$), we can derive the following asymptotic distribution.

\begin{theorem}\label{thm:CLT for raw param}
	Under Assumptions \ref{assump:I(1) process}-\ref{assump:max eigenvalues for covariance}, as $T\to\infty$
	\begin{align*}
		 T^{1/4}(\hat{\alpha}_i-\alpha_{0i})\to_D\mathbf{MN}\left(0,\bar{\omega}_{11}^{\theta,i}\frac{\alpha_{0i}\alpha_{0i}'}{\|\alpha_{0i}\|^2}\right)~\text{and}~ T^{3/4}\left(\hat{\alpha}_i^{\circ}-\frac{\alpha_{0i}}{\|\alpha_{0i}\|}\right)\to_D\mathbf{MN}\left(0,\frac{Q_i^{(2)}\bar{\omega}_{22}^{\theta,i}Q_i^{(2)'}}{\|\alpha_{0i}\|^2}\right).
	\end{align*}
\end{theorem}
Here, the normalization of $\hat{\alpha}_i$ scales it to the unit sphere, focusing on angular convergence rather than magnitude. Consequently, the convergence rate is accelerated due to the differing rates for $\hat{\theta}_i$. This suggests that imposing the constraint $\|\alpha_{0i}\|=1$ on the binary probability allows $\hat{\alpha}_i^{\circ}$ to serve as a more precise estimator of $\alpha_{0i}$.

Estimating binary event probabilities is a crucial aspect of statistical analysis. The following corollary presents the corresponding theoretical result.
\begin{corollary}\label{coro:CLT for binary probabilities}
	Under Assumptions \ref{assump:I(1) process}-\ref{assump:max eigenvalues for covariance}, as $T,N\to\infty$
	\begin{align*}
		 \frac{\underline{C}_{NT,t}\left(\Psi(\hat{z}_{it})-\Psi(z_{0it})\right)}{|\dot{\Psi}(z_{0it})|\sqrt{\frac{\underline{C}_{NT,t}^2}{\sqrt{T}}\bar{\omega}_{11}^{\theta,i}(\alpha_{0i}'g_{0it})^2+\frac{\underline{C}_{NT,t}^2t^{\delta}}{N}\lambda_{0i}'\Omega_{f,t}^{-1}\lambda_{0i}}}\to_D\mathbf{N}(0,1),
		\end{align*}
	where $\underline{C}_{NT,t}=\min\{N^{1/2}t^{-\delta/2},T^{1/4}\}$.
\end{corollary}
Corollary \ref{coro:CLT for binary probabilities} indicates that the convergence rate is $\min\{N^{1/2}t^{-\delta/2},T^{1/4}\}$. If $F_0$ is observable, the convergence rate for $\hat{\alpha}_i$ is $T^{1/4}$; If $A_0$ is observable, the rate for $\hat{f}_t$ is $N^{1/2}t^{-\delta/2}$. Therefore, when estimating both parameters simultaneously, the best rate for $\hat{z}_{it}$ is the minimum of $T^{1/4}$ and $N^{1/2}t^{-\delta/2}$.

To extend the previous results to the plug-in version, we first derive estimators for the quantities of interest. According to Theorems  \ref{thm:asymptotic distribution for theta and f} and \ref{thm:CLT for raw param}, the estimators $\hat{f}_t\sim\mathbf{N}(f_{0t},\frac{t^{\delta}}{N}\Omega_{f,t}^{-1})$ and $\hat{\alpha}_i\sim\mathbf{MN}(\alpha_{0i},\frac{1}{\sqrt{T}}\bar{\omega}_{11}^{\theta,i}\frac{\alpha_{0i}\alpha_{0i}'}{\|\alpha_{0i}\|^2})$. Based on these distributions, we define two estimators for the inverse Hessian of $J_{NT,11}^{(i)}(\alpha_{0i},F_{0})$: $\left[J_{NT,11}^{(i)}(\hat{\alpha}_i,\hat{F})\right]^{-1}$ and $\left[\underline{J}_{NT,11}^{(i)}(\hat{\alpha}_i,\hat{F})\right]^{-1}$. Similarly, for the inverse Hessian of $J_{NT,22}^{(t)}(A_{0},f_{0t})$, we define: $\left[J_{NT,22}^{(t)}(\hat{A},\hat{f}_t)\right]^{-1}$ and $\left[\underline{J}_{NT,22}^{(t)}(\hat{A},\hat{f}_t)\right]^{-1}$. Specifically, we define the dominant terms for the Hessians as:
\begin{align*}
	\underline{J}_{NT,11}^{(i)}(\hat{\alpha}_i,\hat{F})=-\sum_{t=1}^TK(\hat{z}_{it})\hat{g}_{it}\hat{g}_{it}'\quad\text{and}\quad \underline{J}_{NT,22}^{(t)}(\hat{A},\hat{f}_t)=-\sum_{i=1}^NK(\hat{z}_{it})\hat{\lambda}_{i}\hat{\lambda}_{i}',
\end{align*}
where $\hat{z}_{it}=\hat{\beta}_{i}'x_{it}+\hat{\lambda}_i'\hat{g}_t$ and $\hat{g}_t=(x_{it}',\hat{f}_t')'$. The following corollary establishes the consistency of the proposed estimators:
\begin{corollary}\label{coro:consistency for plug-in estimators}
	Under Assumptions \ref{assump:I(1) process}-\ref{assump:max eigenvalues for covariance}, as $T\to\infty$, the following results hold:
	\begin{align*}
		&-\sqrt{T}\left[J_{NT,11}^{(i)}(\hat{\alpha}_i,\hat{F})\right]^{-1},\quad-\sqrt{T}\left[\underline{J}_{NT,11}^{(i)}(\hat{\alpha}_i,\hat{F})\right]^{-1}\to_P\bar{\omega}_{11}^{\theta,i}\frac{\alpha_{0i}\alpha_{0i}'}{\|\alpha_{0i}\|^2},\\
		&-Nt^{-\delta}\left[J_{NT,22}^{(t)}(\hat{A},\hat{f}_t)\right]^{-1},\quad -Nt^{-\delta}\left[\underline{J}_{NT,22}^{(t)}(\hat{A},\hat{f}_t)\right]^{-1}\to_P\Omega_{f,t}^{-1}.
	\end{align*}
\end{corollary}
Based on Corollary \ref{coro:consistency for plug-in estimators} and Slutsky's theorem, we replace the limiting distribution of $\hat{\alpha}_i$ in Theorem \ref{thm:CLT for raw param} and that of $\hat{f}_t$ in Theorem \ref{thm:asymptotic distribution for theta and f} to obtain the following plug-in version of the limiting distribution:
\begin{align*}
	\left[-\underline{J}_{NT,11}^{(i)}(\hat{\alpha}_i,\hat{F})\right]^{1/2}(\hat{\alpha}_i-\alpha_{0i})\to_D\mathbf{N}\left(0,I_{q+r}\right)~\text{and}~ \left[-\underline{J}_{NT,22}^{(t)}(\hat{A},\hat{f}_t)\right]^{1/2}(\hat{f}_t-f_{0t})\to_D\mathbf{N}(0,I_r).
\end{align*}

Next, we provide the local time estimator for $L_{1i}(1,0)$:
\begin{align}\label{eq:plug-in local time}
	\hat{L}_{1i}(1,0)=\frac{\|\hat{\alpha}_i\|}{\sqrt{T}}\sum_{t=1}^T\dot{\Psi}(\hat{z}_{it}).
\end{align}
This estimator is intuitive because, by consistency, $\frac{\|\alpha_{0i}\|}{\sqrt{T}}\sum_{t=1}^T\Psi'(z_{0it})\to_DL_{1i}(1,0)\Psi(z)|_{-\infty}^{\infty}=L_{1i}(1,0)$. Finally, we define the mean squared error (MSE) for the observations to assess the accuracy of the model estimates:
\begin{align*}
	\hat{\text{MSE}}=\frac{1}{N\sqrt{T}}\sum_{i=1}^N\sum_{t=1}^T\left[y_{it}-\Psi(\hat{z}_{it})\right]^2.
\end{align*}

The following proposition establishes the consistency of the proposed estimators:
\begin{proposition}\label{prop:consistency for MSE}
	Under Assumptions \ref{assump:I(1) process}-\ref{assump:max eigenvalues for covariance}, as $T\to\infty$, the following results hold:
	\begin{align*}
        \hat{L}_{1i}(1,0)\to_PL_{1i}(1,0)\quad\text{and}\quad \hat{\text{MSE}}\to_P\int_{\mathbb{R}}\Psi(s)(1-\Psi(s))ds\left(\mathrm{plim}_{N\to\infty}\frac{1}{N}\sum_{i=1}^N\frac{L_{1i}(1,0)}{\|\alpha_{0i}\|}\right).
	\end{align*}
\end{proposition}
The local time estimator in Equation \eqref{eq:plug-in local time} is not unique; it suffices that the nonlinear function within it is integrable and integrates to 1 over $\mathbb{R}$. Additionally, $\hat{\text{MSE}}$ serves as an approximation of $\frac{1}{N\sqrt{T}}\log L(\hat{B},\hat{\Lambda},\hat{F})$ in \eqref{eq:log-likelihood function} (see \citealt{gao2023binary} for this insight) assessing the model's goodness of fit.

\begin{remark}
	When allowing for partial cointegration in $z_{0it}$---where the cointegration rank is smaller than $q+r-1$---in the series $\{g_{0it}\}$, while maintaining the nonstationarity of $\{z_{0it}\}$, similar results can be obtained. For instance, if $\lambda_i'f_{0t}$ is stationary for all $i=1,...,N$, the primary characteristics of the model remain largely unchanged, except for $H_{1i}$ and $H_{2i}$, since $\{z_{0it}\}$ continuous to be nonstationary. In this scenario, the local time corresponds to the Brownian motion $\beta_{0i}'E_i(s)/\|\alpha_{0i}\|$, rather than $(\beta_{0i}'E_i(s)+\lambda_{0i}'V(t))/\|\alpha_{0i}\|$, leading to a modification in $H_{2i}$. To uphold the asymptotic theory under these conditions, it is essential that $\{z_{0it}\}$ remains nonstationary and that the cointegration rank of each $\{g_{0it}\}$ is less than $q+r-1$, thereby ensuring the persistence of dual convergence rates.
\end{remark}

\subsection{Selecting the Number of Factors}\label{subsec:Selecting the Number of Factors}
In this section we introduce a rank minimization method to select the number of factors.

Assume $k$ is a positive integer larger than $r$. We solve the optimization problem  in Equation \eqref{eq:MLE estimators} using $k$ factors, resulting in estimators $\hat{B}^k$, $\hat{\Lambda}^k$ and $\hat{F}^k$. Define $(\hat{\Lambda}^k)'\hat{\Lambda}^k/N=\text{diag}(\hat{\sigma}_{N,1}^k,...,\hat{\sigma}_{N,k}^k)$. The estimator for the number of factors is then given by
\begin{align*}
	\hat{r}=\sum_{j=1}^k\mathbb{I}_{\{\hat{\sigma}_{N,j}^k>\pi_{NT}\}},
\end{align*}
where $\pi_{NT}$ is a sequence approaching zero as $N,T\to\infty$. In other words, $\hat{r}$ counts the number of diagonal elements in $(\hat{\Lambda}^k)'\hat{\Lambda}^k/N$ that exceed the threshold $\pi_{NT}$. To elucidate, decompose $\hat{\Lambda}^k$ into $\hat{\Lambda}^k=(\hat{\Lambda}^{k,r},\hat{\Lambda}^{k,-r})$, where $\hat{\Lambda}^{k,r}$ comprises the first $r$ columns of $\hat{\Lambda}^k$, and $\hat{\Lambda}^{k,-r}$ includes the remaining $k-r$ columns. It can be shown that $\hat{\sigma}_{N,j}^k=\sigma_{N,j}^k+o_P(1)\to_P\sigma_j$ for $j=1,...,r$, and $\hat{\sigma}_{N,j}^k=o_P(1)$ for $j>r$. Consequently, $(\hat{\Lambda}^k)'\hat{\Lambda}^k/N$ converges in probability to a matrix with rank $r$ at certain rates. By selecting a suitable threshold $\pi_{NT}$, which is greater than this rate and less than $\sigma_r$, we can accurately determine $r$. 
\begin{theorem}\label{thm:number of factors}
	Under Assumptions \ref{assump:I(1) process}-\ref{assump:max eigenvalues for covariance}, as $N,T\to\infty$, if $k>r$, $\sqrt{T}/N=o(1)$, $\pi_{NT}\to0$, and $\pi_{NT}C_{NT}^2T^{-1/2}\to\infty$, then $P(\hat{r}=r)\to1$.
\end{theorem}
A threshold value of
$
	\pi_{NT}=\hat{\sigma}_{N,1}^k\left(C_{NT}^2T^{-1/2}\right)^{-1/3}
$
has been found to perform well in numerical simulations. Some methods of choosing the number of factors with the help of eigenvalue properties can also be generalized here, e.g. \cite{trapani2018randomized} and \cite{yu2024testing}.

\subsection{Cointegrated Single Index}\label{subsec:Cointegrated Single-Indexes}
In this subsection, we examine the case where a linear cointegration happens among components of $g_{0it}$ for all $i=1,...,N$. In other words $\alpha_{0i}'g_{0it}\sim I(0)$ for every $i$. 

Within this context, we solve the optimization problem in Equation \eqref{eq:MLE estimators} under the assumption that $g_{0it}$ is a $(q+r)$-dimensional $I(1)$ process and that the single index $\alpha_{0i}'g_{0it}$ is $I(0)$. In Remark \ref{rmk:remark 4} below, we relax the assumption to accommodate nonstationarity in the single indices $\{z_{0it}\}$ for some $i$.

To derive the asymptotic properties of the estimators $(\hat{A},\hat{F})$, we first rotate the coordinate system using a $(q+r)\times(q+r)$ orthogonal matrix $Q_i=(Q_i^{(1)},Q_i^{(2)})$, where $Q_i^{(1)}=\alpha_{0i}/\|\alpha_{0i}\|$ defines the primary axis. This transformation enables us to express the single index $z_{0it}$ as
\begin{align*}
	 z_{0it}=\alpha_{0i}'Q_iQ_i'g_{0it}=\theta_{0i}^{(1)}h_{0it}^{(1)}+\theta_{0i}^{(2)'}h_{0it}^{(2)},
\end{align*}
where $\theta_{0i}^{(1)}=\|\alpha_{0i}\|$ and $\theta_{0i}^{(2)}=0$. In contrast to Section \ref{subsec:Theoretical Results for Nonstationary Models}, here the component $h_{0it}^{(1)}=\alpha_{0i}'g_{0it}/\|\alpha_{0i}\|$ is a stationary scalar process, while $h_{0it}^{(2)}=Q_i^{(2)'}g_{0it}$ is a $(q+r-1)$-dimensional nonstationary process. Because the series $\{z_{0it}\}$ is concentrated within the effective range of $K$, some aspects of the original asymptotic theory must be revised, and we update the assumptions.

\begin{assumption}\label{assump:cointegrated single-indexes}(Cointegrated Single Indices)
	\begin{enumerate}[(i)]
		\item There exists a set $\Xi=[\Xi_l,\Xi_u]$ such that $z_{0it}$ belongs to $\Xi$ w.p.a.1. Moreover, $\Psi(\Xi_l)>0$ and $\Psi(\Xi_u)<1$.
		\item $K$, $\dot{\Psi}$, $\dot{M}\dot{\Psi}$, $M\ddot{\Psi}$, $\ddot{M}\Psi^{1/2}(1-\Psi)^{1/2}$, $\dot{M}\Psi^{1/2}(1-\Psi)^{1/2}$, and $M^3\dot{\Psi}$ all belong to the class $\mathbb{F}_B$.
	    \item $\varPi^e(1)$ and $\varPi^v(1)$ are nonsingular.  $\varPi(1)$ has rank $q+r-1$ and $\alpha_{0i}'\varPi(1)=0_{1\times (q+r)}$, where $\varPi(L)=\sum_{k=0}^{\infty}\text{diag}(\varPi^e_k,\varPi^v_k)L^k$. All other conditions remain as in Assumption \ref{assump:I(1) process}.
	    \item Let $\{z_{0it}\}$ be a strictly stationary process that is $\alpha$-mixing over $t$ with mixing coefficient $\alpha_{ij}(\tau)$ satisfying $\max_{i\geq1}\sum_{\tau=1}^{\infty}(\alpha_{ii}(\tau))^{\nu/(4+\nu)}<\infty$, $\frac{1}{N}\sum_{i,j=1}^N\sum_{\tau=1}^{\infty}(\alpha_{ij}(\tau))^{\nu/(4+\nu)}<\infty$, and $\frac{1}{N}\sum_{i,j=1}^N(\alpha_{ij}(0))^{\nu/(4+\nu)}<\infty$ for some $\nu>0$.
	    \item Assume that $\max_{i\geq,t\geq 1}E\|g_{0it}\|^{4+\nu}<\infty$ and $\|g_{0it}\|\leq C$ for some positive constant $C$.
	    \item For each $t$, as $N\to\infty$, $\frac{1}{\sqrt{N}}\sum_{i=1}^NM(z_{0it})\lambda_{0i}u_{it}\to_D\mathbf{N}(0,\Omega_{f,t}^*)$, where the covariance matrix $\Omega_{f,t}^*=\mathrm{lim}_{N\to\infty}\frac{1}{N}\sum_{i=1}^NE[K(z_{0it})]\lambda_{0i}\lambda_{0i}'$.
	\end{enumerate}
\end{assumption}
Assumption \ref{assump:cointegrated single-indexes}(i) ensures that the support of $\{z_{0it}\}$ is bounded, which is reasonable when $z_{0it}$ is a stationary scalar.  Assumption \ref{assump:cointegrated single-indexes}(ii) relaxes the function class in Assumption \ref{assump:some integrable functions}. Assumption \ref{assump:cointegrated single-indexes}(iii) follows from the requirement only $\{z_{0it}\}$ is stationary. Assumptions \ref{assump:cointegrated single-indexes} (iv) and (v) ensure the $\alpha$-mixing for $z_{0it}$ and the boundedness of $g_{0it}$, such as \cite{trapani2021inferential}.

We first give results for the average rate of convergence and the number of factors under the cointegrated single-index case. Define $\underline{D}_T=\text{diag}(\sqrt{T},TI_{q+r-1})$. The estimation of the number of factors is fully consistent with Section 3.2 except for the choice of thresholds $\pi_{NT}$.

\begin{theorem}\label{thm:average convergence rate and number for cointegrated single-indexes}
	Under Assumptions \ref{assump:Factors and loadings}, \ref{assump:max eigenvalues for covariance}, and \ref{assump:cointegrated single-indexes}, the following results hold. \begin{enumerate}[(i)]
		\item $\frac{1}{\sqrt{N}}\|(\hat{\Theta}-\Theta_0)\underline{D}_T\|=O_P\left(C_{NT}^{-1}\right)\quad\text{and}\quad\frac{1}{\sqrt{T}}\|\hat{F}-F_0\|=O_P\left(C_{NT}^{-1}\right).$
		\item If $k>r$, $\pi_{NT}\to0$, and $\pi_{NT}C_{NT}^2\to\infty$, then $P(\hat{r}=r)\to1$.
	\end{enumerate}
\end{theorem}
We find that the convergence rates of both the coefficients and the factors are improved, suggesting that estimating the model in the cointegrated single-index case is more accurate. In addition, the threshold selection for the number of estimated factors is less demanding.

We now study the asymptotic distribution of the estimators $(\hat{\theta}_i,\hat{f}_t)$. Unlike in Section \ref{subsec:Theoretical Results for Nonstationary Models}, where $h_{0it}^{(1)}$ is a nonstationary process, here $h_{0it}^{(1)}$ is stationary. This change may affect the convergence rate of $\hat{\theta}_{i}^{(1)}$, which lies in the direction of $\alpha_{0i}$. Similarly, the convergence rate of $\hat{\theta}_{i}^{(2)}$ in the orthogonal direction may also differ. The following theorem addresses these questions.
\begin{theorem}\label{thm:asymptotic theory for cointegrated single-indexes}
	Under Assumptions \ref{assump:Factors and loadings}, \ref{assump:max eigenvalues for covariance}, and \ref{assump:cointegrated single-indexes}, as $N,T\to\infty$,
	\begin{align*}
		 \underline{D}_T(\hat{\theta}_i-\theta_{0i})\to_D\Omega_{\theta,i}^{*-1}\xi_{i},\quad\text{and}\quad\sqrt{N}(\hat{f}_t-f_{0t})\to_D\mathbf{N}(0,\Omega_{f,t}^{*-1}),
	\end{align*}
where $\xi_{i}=(\xi_{1i},\xi_{2i}')'$ with
\begin{align*}
	 \xi_{1i}=&\sqrt{E\left[M(\|\alpha_{0i}\|h_{0i1}^{(1)})h_{0i1}^{(1)}\right]^2}\int_{0}^1dU_i(s)~\text{and}~\xi_{2i}=\sqrt{E\left[M(\|\alpha_{0i}\|h_{0i1}^{(1)})\right]^2}\int_{0}^1H_{2i}(s)dU_i(s),\\
	\Omega_{\theta,i}^*=&\begin{pmatrix}
		 E\left[K(\|\alpha_{0i}\|h_{0i1}^{(1)})(h_{0i1}^{(1)})^2\right]&E\left[K(\|\alpha_{0i}\|h_{0i1}^{(1)})h_{0i1}^{(1)}\right]\int_0^1H_{2i}'(s)ds\\
		 E\left[K(\|\alpha_{0i}\|h_{0i1}^{(1)})h_{0i1}^{(1)}\right]\int_0^1H_{2i}(s)ds&E\left[K(\|\alpha_{0i}\|h_{0i1}^{(1)})\right]\int_0^1H_{2i}(s)H_{2i}'(s)ds
	\end{pmatrix},
\end{align*}
and $\Omega_{f,t}^*=\lim_{N\to\infty}\frac{1}{N}\sum_{i=1}^NE\left(K\left(\|\alpha_{0i}\|h_{0it}^{(1)}\right)\right)\lambda_{0i}\lambda_{0i}'$ and $H_{2i}=Q_i^{(2)'}H_i(r)$.
\end{theorem}
Theorem \ref{thm:asymptotic theory for cointegrated single-indexes} shows that the convergence rates for $\hat{\theta}_i^{(1)}$ and $\hat{\theta}_i^{(2)}$ differ from those in Theorem \ref{thm:asymptotic distribution for theta and f}. Specifically, when the single indices are cointegrated, the convergence rates of the parameter estimators improve, and the asymptotic results resemble those observed in linear models. For $\hat{f}_t$, the conventional asymptotics hold, owing to a constant lower bound on $K(z_{0it})$ for all $i=1,...,N$ and $t=1,...,T$.

Rewrite $\Omega_{\theta,i}^*=\begin{pmatrix} \omega_{11}^{\theta,i*}& \omega_{12}^{\theta,i*}\\\omega_{21}^{\theta,i*}&\omega_{22}^{\theta,i*}	 \end{pmatrix}$ and the inverse matrix $\Omega_{\theta,i}^{*-1}=\begin{pmatrix} \bar{\omega}_{11}^{\theta,i*}& \bar{\omega}_{12}^{\theta,i*}\\\bar{\omega}_{21}^{\theta,i*}&\bar{\omega}_{22}^{\theta,i*}	 \end{pmatrix}$, where
\begin{align*}
	 \bar{\omega}_{11}^{\theta,i*}=\left(\omega_{11}^{\theta,i*}-\omega_{12}^{\theta,i*}(\omega_{22}^{\theta,i*})^{-1}\omega_{21}^{\theta,i*}\right)^{-1}\quad\text{and}\quad\bar{\omega}_{22}^{\theta,i*}=\left(\omega_{22}^{\theta,i*}-\omega_{21}^{\theta,i*}(\omega_{11}^{\theta,i*})^{-1}\omega_{12}^{\theta,i*}\right)^{-1}.
\end{align*}

Leveraging the linear relationship between $\hat{\alpha}_i$ and $\hat{\theta}_i$, we can derive the asymptotic results for $\hat{\alpha}_i$. The following theorem presents the asymptotic distributions of both $\hat{\alpha}_i$ and $\hat{\alpha}_i^{\circ}$.

\begin{theorem}\label{thm:asymptotic theory for cointegrated single-indexes in alpha}
	Under Assumptions \ref{assump:Factors and loadings}, \ref{assump:max eigenvalues for covariance}, and \ref{assump:cointegrated single-indexes}, as $N,T\to\infty$,
		\begin{align*}
		 T^{1/2}(\hat{\alpha}_i-\alpha_{0i})\to_D&\frac{\alpha_{0i}}{\|\alpha_{0i}\|}\left[\bar{\omega}_{11}^{\theta,i*}\xi_{1i}+\bar{\omega}_{12}^{\theta,i*}\xi_{2i}\right],~T\left(\hat{\alpha}_i^{\circ}-\frac{\alpha_{0i}}{\|\alpha_{0i}\|}\right)\to_D\frac{Q_i^{(2)}}{\|\alpha_{0i}\|}\left[\bar{\omega}_{21}^{\theta,i*}\xi_{1i}
		+\bar{\omega}_{22}^{\theta,i*}\xi_{2i}\right].
	\end{align*}
where $\xi_{1i}$ and $\xi_{2i}$ are defined in Theorem \ref{thm:asymptotic theory for cointegrated single-indexes}.
\end{theorem}
The convergence rates of the estimators $\hat{\alpha}_i$ and $\hat{\alpha}_i^{\circ}$, as presented in Theorem \ref{thm:asymptotic theory for cointegrated single-indexes in alpha}, are influenced by the convergence rates of $\hat{\theta}_i^{(1)}$ and $\hat{\theta}_i^{(2)}$, respectively, as detailed in Theorem \ref{thm:asymptotic theory for cointegrated single-indexes}. Notably, these convergence rates for $\hat{\alpha}_i$ and $\hat{\alpha}_i^{\circ}$ differ markedly from those in Theorem \ref{thm:CLT for raw param}, primarily due to the impact of cointegration. Consequently, the asymptotic distributions undergo significant alterations.

\begin{remark}\label{rmk:remark 4}
	In practice, the time series $\{z_{0it}\}$ may have both stationary and nonstationary series for different $i$. We can partition the indices into two exclusive sets, $\mathbb{N}_1$ and $\mathbb{N}_2$, such that ${z_{0it}}$ is stationary for $i \in \mathbb{N}_1$ and nonstationary for $i \in \mathbb{N}_2$. In this scenario, it is necessary to integrate the asymptotic theories in Sections \ref{subsec:Theoretical Results for Nonstationary Models} and \ref{subsec:Cointegrated Single-Indexes}.
	
	For the estimators $\hat{\theta}_i$ and $\hat{\alpha}_i$, the asymptotic distributions remain consistent within their separate sets. However, the estimation of $\hat{f}_t$ becomes more intricate. When $t$ is large, the asymptotic behavior of $\hat{f}_t$ is predominantly influenced by the stationary components in $\mathbb{N}_1$. While for moderate values of $t$, the asymptotic properties are determined by the combined contributions of both $\mathbb{N}_1$ and $\mathbb{N}_2$. Identifying $\mathbb{N}_1$ and $\mathbb{N}_2$ is not easy, we leave it to our future research works. 
\end{remark}

\section{Simulation}\label{sec:Simulation}
In this section, we perform simulations to verify the accuracy of the estimation results. 

\subsection{Simulation Design}\label{subsec:Simulation Design}
To do so, we consider a model with the following data generating process (DGP).

\begin{enumerate}[\textbf{Case} 1.]
	\item \textbf{Non-stationary probabilities}: $q=4$, $r=2$, $x_{it}=x_{it-1}+e_{it}$ with $e_{it}=0.1\times e_{it-1}+0.1\times \mathbf{N}(0,I_4)$, $f_{0t}=f_{0t-1}+v_{t}$ with $v_t=0.01\times \mathbf{N}(0,I_2)$, $\beta_{0i}\sim U[0,1]$, and $\lambda_{0i}\sim\mathbf{N}(0,\mathrm{diag}(2,1))$. The error term $\epsilon_{it}$ is linked to the binary response via either the logit or probit function. The covariate $\{x_{it}\}$ is observable, while the remaining parameters are unobservable.
	\item \textbf{Cointegrated probabilities}: $q=4$, $r=2$, $f_{0t}=f_{0t-1}+v_{t}$ with $v_t=0.01\times \mathbf{N}(0,I_2)$, $\beta_{0i}=(1,0.5,0.5,1)'$, $\lambda_{0i}\sim\mathbf{N}(0,\mathrm{diag}(2,1))$, and $x_{it}=(-0.5\lambda_{0i}(1)\times f_{0t}(1)-0.25\times \lambda_{0i}(2)f_{0t}(2),-0.5\lambda_{0i}(1)\times f_{0t}(1),-0.5\times \lambda_{0i}(2)f_{0t}(2),-0.25\lambda_{0i}(1)\times f_{0t}(1)-0.5\times \lambda_{0i}(2)f_{0t}(2))'+e_{it}$ with $e_{it}=0.1\times e_{it-1}+\mathbf{N}(0,I_4)$, where $\lambda_{0i}(j)$ and $f_{0t}(j)$ represent the $j$th value of $\lambda_{0i}$ and $f_{0t}$, respectively. As in Case 1, the binary response is modeled using either the logit or probit function. The covariate $\{x_{it}\}$ is observable, while all other parameters remain unobservable.
\end{enumerate}

For each generated dataset, we first determine the number of factors. Then, we evaluate the parameter estimates by measuring the error according to the following criteria, where $M$ denotes the number of iterations.
\begin{align}\label{eq:MAE}
	\begin{aligned}
	\text{MAE 1} = &\frac{1}{NTM}\sum_{i=1}^N\sum_{t=1}^T\sum_{j=1}^M|\hat{z}_{it}^{(j)}-z_{0it}|,~~\text{MAE 2} = \frac{1}{NTM}\sum_{i=1}^N\sum_{t=1}^T\sum_{j=1}^M|\hat{\beta}_{i}^{(j)'}x_{it}-\beta_{0i}'x_{it}|,\\
	\text{MAE 3} = &\frac{1}{NTM}\sum_{i=1}^N\sum_{t=1}^T\sum_{j=1}^M|\hat{\lambda}_{i}^{(j)'}\hat{f}_{t}^{(j)}-\lambda_{0i}'f_{0t}|,~~\text{MAE 4} = \frac{1}{NTM}\sum_{i=1}^N\sum_{t=1}^T\sum_{j=1}^M\|\hat{\beta}_{i}^{(j)}-\beta_{0i}\|,
	\end{aligned}
\end{align}
 where $X^{(j)}$ represents the $j$th replication for $X\in\{(z_{0it},\beta_{0i},\lambda_{0i},f_{0t})\}$. We set the number of simulations $M$ to 200.

\subsection{Simulation Results}\label{subsec:Simulation Results}

\begin{table}[htb]\footnotesize
	\setlength\tabcolsep{2pt}
	\centering
	\caption{\textit{Notes}. ``Logit'' and ``Probit'' refer to the logit and probit link functions, respectively, while ``Nonstationary'' and ``Cointegration'' denote the single-index scenarios of being nonstationary and cointegrated.}\label{tab:Numerical simulation results}
	\begin{tabular*}{\textwidth}{@{\extracolsep{\fill}}rccccrcccrcccrccc}
		\toprule
		&       & \multicolumn{7}{c}{Nonstationary}                     &       & \multicolumn{7}{c}{Cointegration} \\
		\cmidrule{3-9}\cmidrule{11-17}          &       & \multicolumn{3}{c}{Logit} &       & \multicolumn{3}{c}{Probit} &       & \multicolumn{3}{c}{Logit} &       & \multicolumn{3}{c}{Probit} \\
		\cmidrule{3-5}\cmidrule{7-9}\cmidrule{11-13}\cmidrule{15-17}          & $N$\textbackslash$T$ & 100   & 300   & 500   &       & 100   & 300   & 500   &       & 100   & 300   & 500   &       & 100   & 300   & 500 \\
		\midrule
		\multicolumn{1}{c}{\multirow{3}[1]{*}{$\hat{r}$}} & 100   & 1.6985 & 1.9426 & 1.8013 &       & 1.2055 & 1.1410 & 1.2543 &       & 1.3652 & 1.4218 & 1.8149 &       & 1.1323 & 1.6886 & 1.7258 \\
		& 300   & 1.7983 & 1.9894 & 1.9878 &       & 1.8142 & 1.9827 & 1.9352 &       & 1.7935 & 1.9209 & 1.9919 &       & 1.5960 & 1.9142 & 1.9827 \\
		& 500   & 1.8030 & 1.9958 & 2.0065 &       & 1.8930 & 1.9948 & 2.0843 &       & 1.8428 & 1.9252 & 2.0013 &       & 1.7025 & 1.924- & 2.0110 \\
		&       &       &       &       &       &       &       &       &       &       &       &       &       &       &       &  \\
		\multicolumn{1}{c}{\multirow{3}[0]{*}{MAE 1}} & 100   & 0.8970 & 0.5798 & 0.6709 &       & 5.3008 & 7.8874 & 3.3996 &       & 2.5923 & 1.7338 & 0.7315 &       & 1.4395 & 0.5834 & 0.3640 \\
		& 300   & 0.6199 & 0.4216 & 0.4242 &       & 0.4841 & 0.3922 & 0.4350 &       & 0.8119 & 0.3842 & 0.3126 &       & 0.6141 & 0.2840 & 0.2278 \\
		& 500   & 0.5752 & 0.4016 & 0.3749 &       & 0.4138 & 0.4053 & 0.3981 &       & 0.6687 & 0.3827 & 0.2943 &       & 0.5216 & 0.2721 & 0.2131 \\
		&       &       &       &       &       &       &       &       &       &       &       &       &       &       &       &  \\
		\multicolumn{1}{r}{\multirow{3}[0]{*}{MAE 2}} & 100   & 0.4975 & 0.3781 & 0.3739 &       & 0.4672 & 0.6724 & 0.7863 &       & 0.6519 & 0.3581 & 0.2545 &       & 0.8636 & 0.4060 & 0.2398 \\
		& 300   & 0.4525 & 0.3511 & 0.4166 &       & 0.4120 & 0.4145 & 0.6068 &       & 0.5178 & 0.2584 & 0.1964 &       & 0.4985 & 0.2033 & 0.1537 \\
		& 500   & 0.4367 & 0.3256 & 0.3347 &       & 0.3765 & 0.4612 & 0.4688 &       & 0.4912 & 0.2524 & 0.1909 &       & 0.4293 & 0.1972 & 0.1461 \\
		&       &       &       &       &       &       &       &       &       &       &       &       &       &       &       &  \\
		\multicolumn{1}{r}{\multirow{3}[0]{*}{MAE 3}} & 100   & 0.7428 & 0.4868 & 0.5682 &       & 5.2045 & 7.8320 & 3.4389 &       & 2.3198 & 1.6066 & 0.6587 &       & 1.0169 & 0.4231 & 0.2851 \\
		& 300   & 0.4488 & 0.3740 & 0.4312 &       & 0.3126 & 0.3646 & 0.5679 &       & 0.5505 & 0.2768 & 0.2366 &       & 0.3416 & 0.1944 & 0.1645 \\
		& 500   & 0.4236 & 0.3187 & 0.3430 &       & 0.3008 & 0.3901 & 0.4119 &       & 0.4247 & 0.2726 & 0.2163 &       & 0.2869 & 0.1763 & 0.1502 \\
		&       &       &       &       &       &       &       &       &       &       &       &       &       &       &       &  \\
		\multicolumn{1}{r}{\multirow{3}[1]{*}{MAE 4}} & 100   & 0.8308 & 0.3090 & 0.2095 &       & 0.6284 & 0.3405 & 0.2460 &       & 0.3611 & 0.1947 & 0.1355 &       & 0.4901 & 0.2273 & 0.1293 \\
		& 300   & 0.7769 & 0.2802 & 0.2064 &       & 0.5800 & 0.2443 & 0.2104 &       & 0.2821 & 0.1367 & 0.1034 &       & 0.2746 & 0.1081 & 0.0813 \\
		& 500   & 0.7459 & 0.2762 & 0.1880 &       & 0.5371 & 0.2526 & 0.1857 &       & 0.2682 & 0.1350 & 0.1010 &       & 0.2361 & 0.1062 & 0.0776 \\
		\bottomrule
	\end{tabular*}
\end{table}
Table \ref{tab:Numerical simulation results} presents the simulation outcomes. We observe that the estimated number of factors is consistent when both $N$ and $T$ are sufficiently large. 

In the nonstationary scenario, when $N$ is relatively small, increasing $T$ can actually lead to larger errors and less stable results for parameter estimation. This finding aligns with Theorem \ref{thm:average convergence rate}, which states that the rate of convergence is $T^{1/4}\left(\min\{\sqrt{N},\sqrt{T}\}\right)^{-1}$; thus, when $N$ is small, a larger $T$ decreases the convergence rate.

In contrast, for the cointegration scenario, larger values of $N$ and $T$ not only reduce the error but also eliminate prior heterogeneity. This result is consistent with Theorem \ref{thm:asymptotic theory for cointegrated single-indexes} and Theorem \ref{thm:asymptotic theory for cointegrated single-indexes in alpha}, which establish a convergence rate of $\left(\min\{\sqrt{N},\sqrt{T}\}\right)^{-1}$.

\section{Empirical Application}\label{sec:Empirical Application}

In this section, we apply our nonstationary binary factor model to high-frequency financial data. By treating daily higher frequency is also possible jumps as binary events and acknowledging that their dynamic is nonstationary (see, for example, \citealt{bollerslev2011estimation} and \citealt{bollerslev2011tails}), we extract the corresponding jump arrival factor and incorporate them into our asset pricing framework:
\begin{align*}
\mathrm{Jump}_{it}=\Psi(\beta_{0i}'x_{it}+\lambda_{0i}'f_{0t})+u_{it},\quad i=1,...,N;\quad t=1,...,T,
\end{align*}
where $\mathrm{Jump}_{it}$ indicates whether or not asset $i$ undergoes jumps on day $t$, with 1 representing presence of jumps and 0 representing absence of jumps.

\subsection{Data}\label{subsec:Data}
We collected intraday observations of S\&P 500 index constituents from January 2004 to December 2016.\footnote{We used the publicly available database provided by \cite{pelger2020understanding}, which includes 332 constituents; see \href{https://doi.org/10.1111/jofi.12898}{https://doi.org/10.1111/jofi.12898}.} Using these high-frequency data, we identify daily jumps with robust detection methods. In particular, we employ the MINRV method proposed by \cite{andersen2012jump}. Additional results regarding the link function and jump detection methods are provided in the Supplementary Material. We set the confidence level for the jump test as 95\%.

Moreover, the overall jump probability (or intensity) trajectory is strongly influenced by volatility (see, for example, \citealt{bollerslev2011tails}). For our covariates, we use the historical volatility of each stock, with daily volatility calculated from high-frequency data.

\subsection{Estimation Results}\label{subsec:Estimation Results}
\begin{figure}[htb]
	\vspace{0pt}\centering
	\includegraphics[scale=0.49]{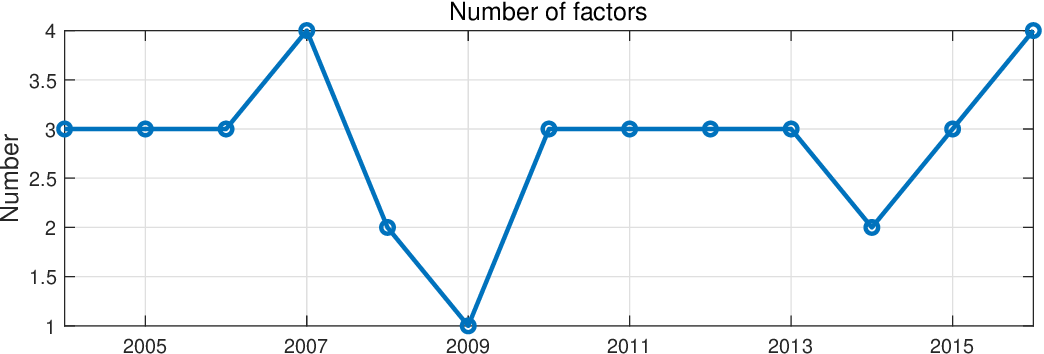}
	\caption{\footnotesize  Estimated number of factors for each year.}\label{fig:Number_Logit_MIN}
\end{figure}

For the period from 2004 to 2016, we estimate three jump arrival factor. To capture dynamic changes, we compute the number of factors for each year, as illustrated in Figure \ref{fig:Number_Logit_MIN}.

Figure \ref{fig:Number_Logit_MIN} shows that the number of factors is typically three, increases by one during the financial crisis, and then drops to one afterward.

Since our jump arrival factor captures the relationship between jump events (such as jump arrivals)—akin to the mutually exciting jumps described in \citealt{dungey2018testing}—it differs significantly from the high-frequency continuous factors in \cite{pelger2020understanding}. While \cite{pelger2020understanding} also constructs jump arrival factors, they rely on sparse jump size data and typically yield only a single factor. In contrast, we extract jump arrival information using a nonlinear factor model that leverages the complete dataset—including both jump occurrence and occurrence rate.

\begin{figure}[htb]
	\vspace{0pt}\centering
	\includegraphics[scale=0.49]{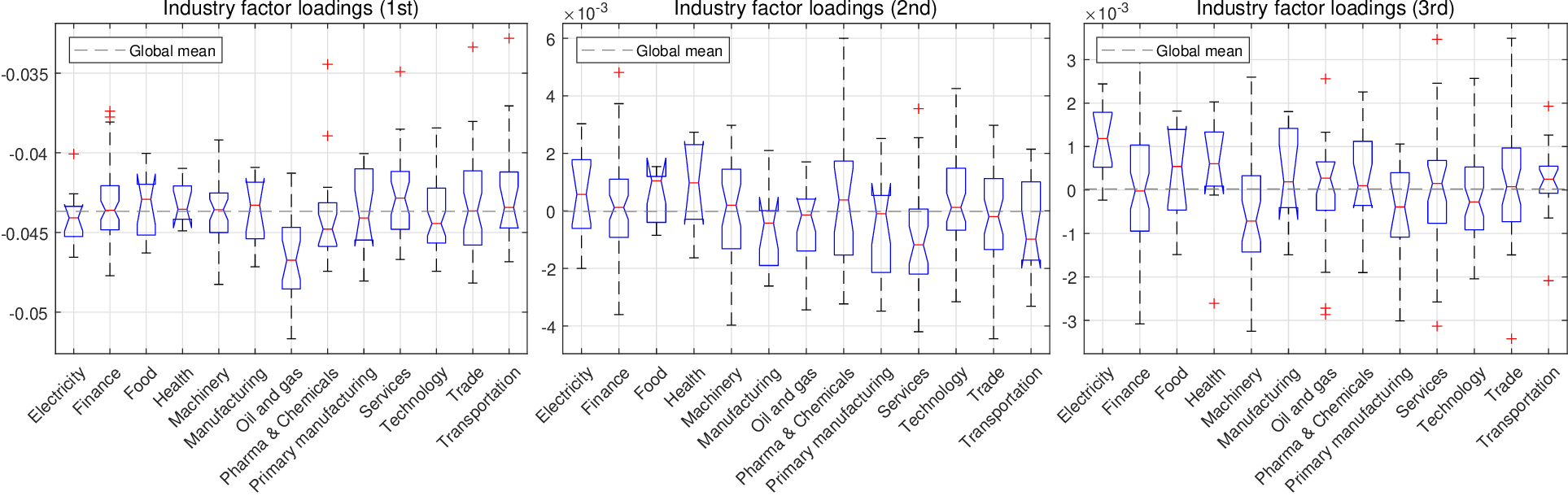}
	\caption{\footnotesize  Box plots of factor loadings across different industries. \textit{Notes}. The three graphs correspond to the three sets of factor loadings, with the horizontal axis representing various industries. The red “+” symbols indicate outliers, and the plots display the confidence interval gaps.}\label{fig:Loadings_Logit_MIN}
\end{figure}
Figure \ref{fig:Loadings_Logit_MIN} shows the distribution of factor loadings across various industries. Our analysis reveals that the first factor's loadings are predominantly negative, with particularly large magnitudes for oil industry. In contrast, the second and third factors fluctuate around zero. Unlike \cite{pelger2020understanding}, which finds that the first four continuous factors are dominated by the financial, oil, and electricity sectors, our results indicate that the industry factor loadings contribute more uniformly to the jump arrival factors.

In addition, we show estimation results for the three jump arrival factors, as well as first-order difference results for the factors.

\begin{figure}[htb]
	\vspace{0pt}\centering
	\includegraphics[scale=0.52]{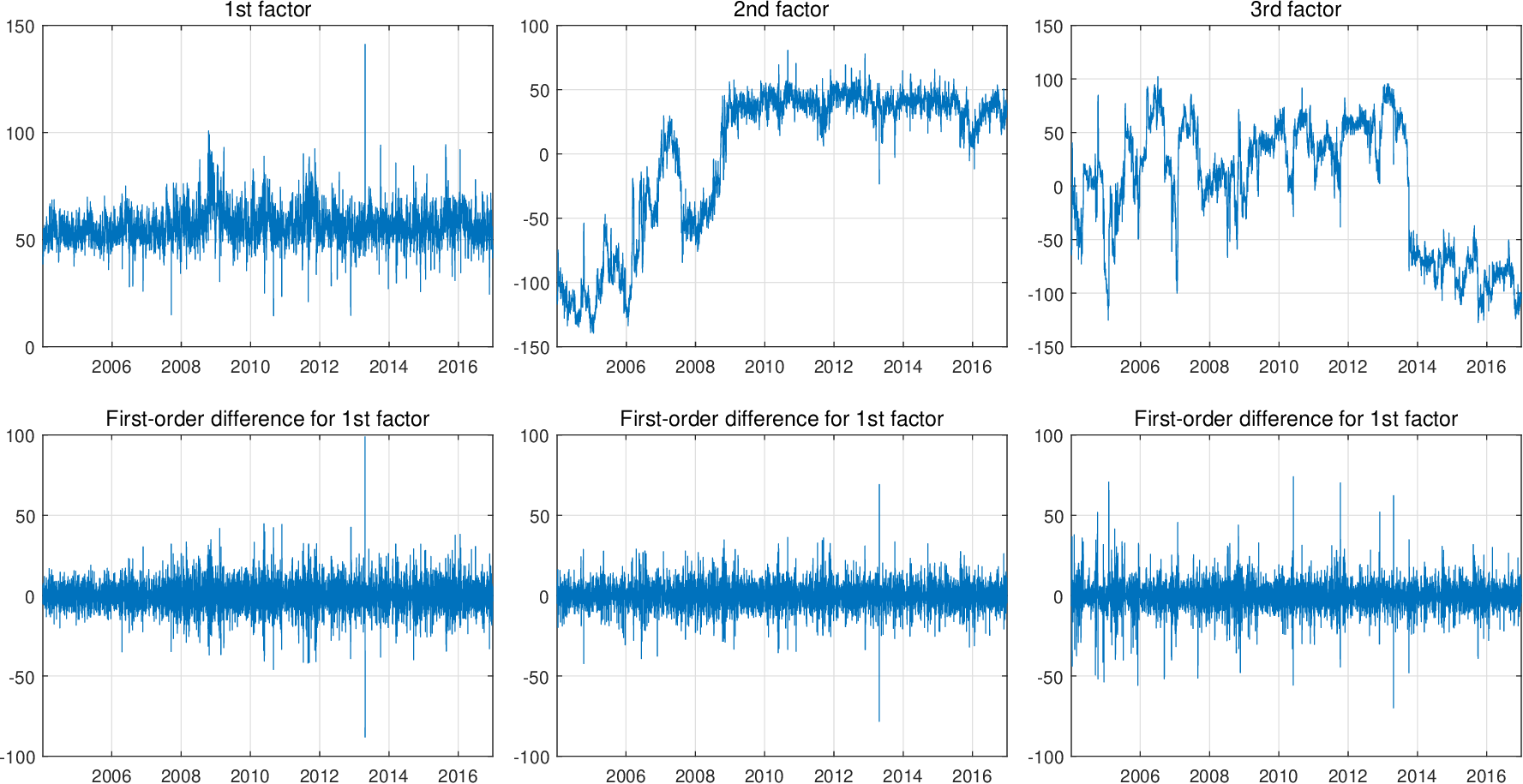}
	\caption{\footnotesize  Estimated factors and their corresponding first-order differences. \textit{Note}. The top panel displays the three estimated factors, while the bottom panel shows the corresponding first-order difference series.}\label{fig:Factors_Logit_MIN}
\end{figure}
The top panel of Figure \ref{fig:Factors_Logit_MIN} displays the three jump arrival factor sequences, which appear to be nonstationary. However, their first-order differences are nearly stationary, aligning well with our model assumptions. To further validate these findings, we conduct an ADF test on the factor series of each year, as presented in Table \ref{tab:ADF results}.

\begin{table}[htb]\footnotesize
	\setlength\tabcolsep{2pt}
	\centering
	\caption{ADF test $p$-value for factors and and their differences. \textit{Notes}. The table's first three rows represent the factors, and the last three rows show their first-order differences.}\label{tab:ADF results}
	\begin{tabular*}{\textwidth}{@{\extracolsep{\fill}}lccccccccccccc}
		\toprule
		Year  & 2004  & 2005  & 2006  & 2007  & 2008  & 2009  & 2010  & 2011  & 2012  & 2013  & 2014  & 2015  & 2016 \\
		\midrule
		1st factor & 0.6002 & 0.4890 & 0.4523 & 0.3966 & 0.4905 & 0.3886 & 0.3605 & 0.4488 & 0.4279 & 0.4736 & 0.4722 & 0.3644 & 0.3734 \\
		2nd factor & 0.6082 & 0.4302 & 0.0527 & 0.3833 & 0.1840 & 0.4706 & 0.4682 & 0.5151 & 0.4812 & 0.2984 & 0.4239 & 0.3021 & 0.5256 \\
		3rd factor & 0.0224 & 0.0252 & 0.2689 & 0.1107 & 0.0017 & 0.4234 & 0.2987 & 0.4244 & 0.5268 & 0.3136 & 0.5261 & 0.5718 & 0.5770 \\
		1st factor diff & 0.0010 & 0.0010 & 0.0010 & 0.0010 & 0.0010 & 0.0010 & 0.0010 & 0.0010 & 0.0010 & 0.0010 & 0.0010 & 0.0010 & 0.0010 \\
		2nd factor  diff & 0.0010 & 0.0010 & 0.0010 & 0.0010 & 0.0010 & 0.0010 & 0.0010 & 0.0010 & 0.0010 & 0.0010 & 0.0010 & 0.0010 & 0.0010 \\
		3rd factor  diff & 0.0010 & 0.0010 & 0.0010 & 0.0010 & 0.0010 & 0.0010 & 0.0010 & 0.0010 & 0.0010 & 0.0010 & 0.0010 & 0.0010 & 0.0010 \\
		\bottomrule
	\end{tabular*}
\end{table}
Table \ref{tab:ADF results} indicates that the first two factors are nonstationary in every year, while the third factor is nonstationary in almost all years. This confirms that the jump arrival factors are predominantly nonstationary. Furthermore, after applying first-order differencing, all factors pass the stationarity test.

Our model has diverse applications. In finance, for instance, we extract jump arrival factors that can be used to screen portfolios, explain asset pricing, and more. The next section demonstrates how our estimated jump arrival factors help explain excess returns, using asset pricing as an example.

\subsection{Applications in Asset Pricing}\label{subsec:Applications in Asset Pricing}
In this section, we investigate the impact of jump arrival factors on pricing models. We analyze each year separately and choose the maximum number of factors observed across all years (i.e., $r=4$) to ensure that variations in the number of factors across periods do not affect the final results.

\begin{figure}[htb]
	\vspace{0pt}\centering
	\includegraphics[scale=0.52]{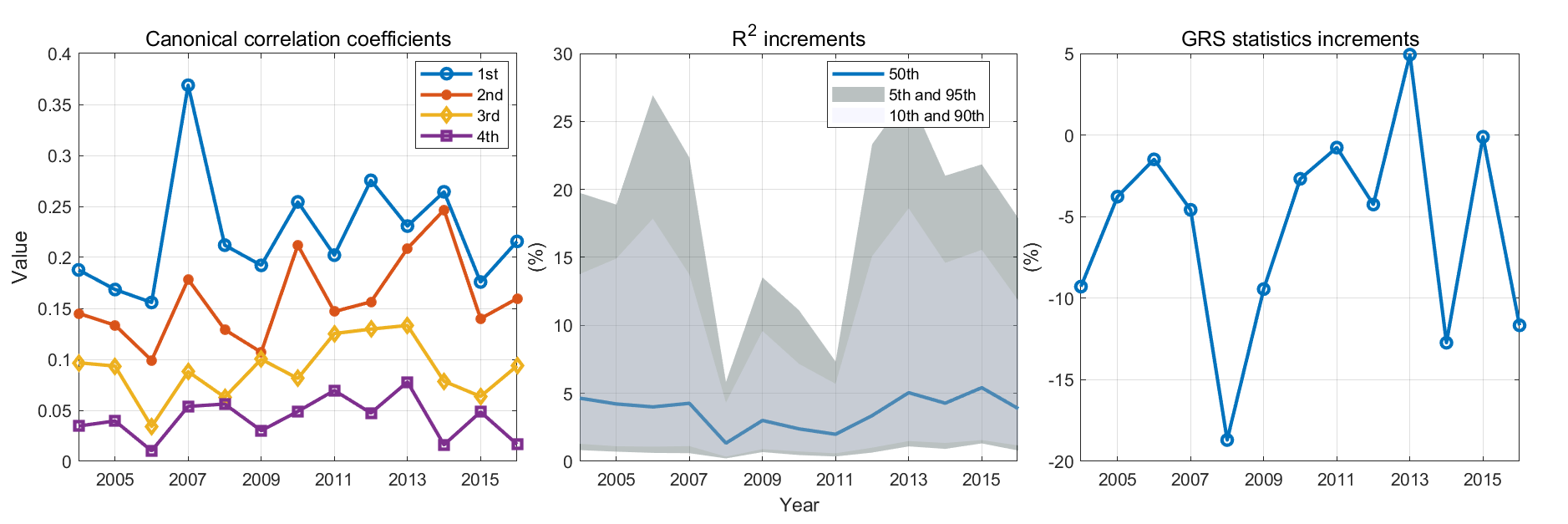}
	\caption{\footnotesize  Canonical correlations and asset pricing results. \textit{Notes}. The left panel displays the canonical correlation coefficients between the four jump arrival factors and the Fama-French-Carhart five factors. The middle panel shows the incremental $R^2$ from adding jump arrival factors to the Fama-French-Carhart five-factor model, while the right panel displays the incremental improvement in the joint test statistic (GRS statistic) for alpha.}\label{fig:APT1_Logit_MIN}
\end{figure}

First, to assess whether the identified jump arrival factors can be explained by established financial factors, we compute the canonical correlations between the four jump arrival factors and the Fama–French–Carhart five factors over the entire sample period. The left panel of Figure \ref{fig:APT1_Logit_MIN} displays these canonical correlation coefficients.

Figure \ref{fig:APT1_Logit_MIN} reveals that most correlation coefficients are relatively low—except for the first coefficient, which shows a modest increase during the financial crisis. This observation suggests that the jump arrival factors are not fully captured by the Fama–French–Carhart five factors.

Motivated by the Fama–French–Carhart five factors model (e.g., \citealt{fama2015five}), we incorporate the jump arrival factors to form the following six-factor model:
\begin{align}\label{eq:FF6}
	\begin{aligned}
	 R_{it}-R_{f,t}=&\alpha_i^*+\beta_{i,MKT}MKT_t+\beta_{i,SMB}SMB_t+\beta_{i,HML}HML_t+\beta_{i,RMW}RMW_t\\
	&+\beta_{i,CMA}CMA_t+\beta_{i,J}'f_t+\epsilon_{i,t}^*,
\end{aligned}
\end{align}
where $R_{it}-R_{f,t}$ represents the excess return of asset $i$ at time $t$ (with $R_{f,t}$ as the risk-free rate), $MKT_t$ is the market excess return, $SMB_t$ captures the size effect, $HML_t$ represents the value factor, $RMW_t$ reflects profitability, $CMA_t$ measures investment, and $f_t$ denotes the set of jump arrival factors.

We evaluate the contribution of the jump arrival factors from two perspectives. First, by comparing the $R^2$ values from regressions based on the six-factor model \eqref{eq:FF6} and the original Fama–French–Carhart five factors model, we determine whether including the jump arrival factors improves the explanation of excess returns. Second, we test the null hypothesis $H_0:\alpha_1^*=...=\alpha_N^*=0$ using the Gibbons–Ross–Shanken (GRS) test to assess the validity of the six-factor model relative to the Fama–French–Carhart five factors model.

The middle panel of Figure \ref{fig:APT1_Logit_MIN} plots the annual increases in $R^2$ for all assets. This panel displays the median, as well as the 5th/95th and 10th/90th percentiles of the $R^2$ increments. Larger $R^2$ improvements indicate that the jump arrival factors enhance the model’s ability to explain asset pricing. On average, the inclusion of jump arrival factors results in nearly a $5\%$ improvement in $R^2$, with some assets exhibiting gains of more than $20\%$.

Similarly, the right panel of Figure \ref{fig:APT1_Logit_MIN} presents the annual changes in the GRS statistic. A lower GRS statistic suggests that the joint alphas are statistically indistinguishable from zero, implying that the model effectively explains excess returns. Negative increments in the GRS statistic indicate that the jump arrival factors enhance the model’s performance. As observed, the GRS statistic decreased in almost every year, reinforcing the effectiveness of the jump arrival factors in explaining asset returns.

\begin{figure}[htbp]
	\vspace{0pt}\centering
	\includegraphics[scale=0.6]{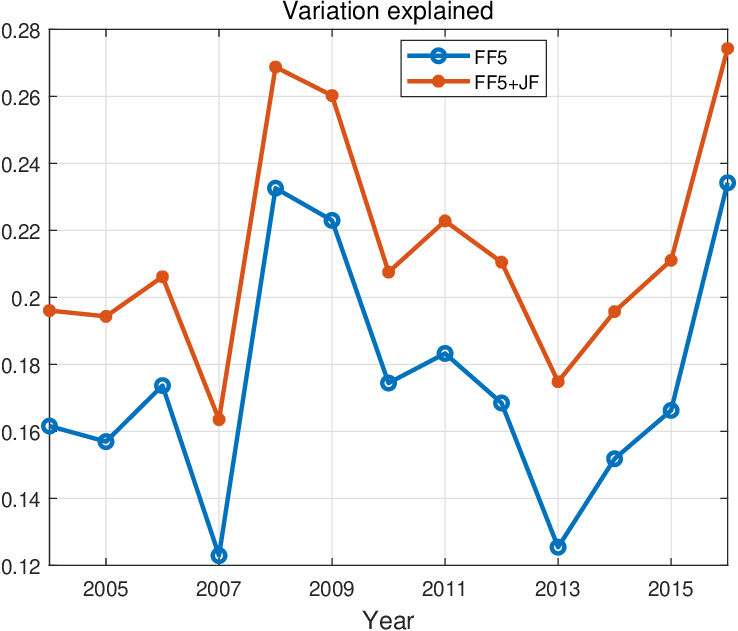}
	\caption{\footnotesize  Time-variation in the percentage of explained variation for different factors. \textit{Note}. This figure plots the percentage of explained variation calculated on a moving window of one year (252 trading days).}\label{fig:APT2_Logit_MIN}
\end{figure}

To further demonstrate that the jump arrival factors have incremental explanatory power, we examine how the proportion of variation explained by jump arrival factors varies over time. We adopt the two-stage regression framework of \cite{fama1973risk}. Figure \ref{fig:APT2_Logit_MIN} illustrates the temporal dynamics using local-regression analysis over a rolling one-year window. The addition of jump arrival factors to the Fama–French–Carhart five factors significantly increases the explained variation by nearly 30\%.

\section{Conclusion}\label{sec:Conclusion}
This paper considers a single-index general factor model with integrated covariates and factors, considering two distinct cases: nonstationary and cointegrated single indices. The estimators are obtained via maximum likelihood estimation, and new asymptotic properties have been established. First, the convergence rates differ between the two cases, with an elevated rate when the single index is cointegrated. Second, while the convergence rate for factor estimators depends on time $t$ in the nonstationary case---necessitating a larger sample size $N$---but is independent of $t$ in the cointegrated case. Third, in a transformed coordinate system, the coefficient estimates exhibit two distinct convergence components. Finally, the limiting distributions of the coefficient estimates are entirely different across the two single-index scenarios. Monte Carlo simulations validate these theoretical results, and empirical studies demonstrate that the extracted nonstationary jump arrival factors play a crucial role in asset pricing. Future research could extend our modelling framework to matrix factor structures, such as \cite{yuan2023two}, \cite{he2024generalized}, and \cite{xu2025quasi}, or to the high-frequency econometrics, such as \cite{pelger2019large} and \cite{chen2024realized}.


\section*{Supplementary Material}
The Supplementary Material contains the proofs of the main theoretical results, additional numerical studies, and more details in the empirical analysis.

\normalem
\bibliographystyle{chicago}

\bibliography{Manuscript}

\newpage
\setcounter{page}{1}

\appendix
\numberwithin{equation}{section}
\numberwithin{theorem}{section}
\numberwithin{lemma}{section}
\renewcommand{\appendixname}{Appendix~\Alph{section}}
\renewcommand\thesection{\Alph{section}}
\renewcommand\theequation{\Alph{section}.\arabic{equation}}



\end{document}